\documentclass[hidelinks,onefignum,onetabnum]{siamart220329}

\usepackage{lipsum}
\usepackage{amsfonts}
\usepackage{graphicx}
\usepackage{epstopdf}
\usepackage{algorithmic}
\usepackage{mathrsfs}
\ifpdf
  \DeclareGraphicsExtensions{.eps,.pdf,.png,.jpg}
\else
  \DeclareGraphicsExtensions{.eps}
\fi

\newsiamremark{remark}{Remark}
\newsiamremark{hypothesis}{Hypothesis}
\crefname{hypothesis}{Hypothesis}{Hypotheses}
\newsiamthm{claim}{Claim}
\newsiamremark{prop}{Proposition}

\headers{}{}

\title{
Iterative Reconstruction Methods for Cosmological X-Ray Tomography}

\author{Julianne Chung\thanks{Department of Mathematics, Emory University, Atlanta, GA 
  (\email{julianne.mei-lynn.chung@emory.edu}), (\email{lonisk@emory.edu}), (\email{yiran.wang@emory.edu}).}
\and Lucas Onisk\footnotemark[1]
\and Yiran Wang\footnotemark[1]}

\def \id {\operatorname{Id}}
\def \comp {\operatorname{comp}}
\def \loc {\text{loc}}

\def \re {\operatorname{Re}}

\def \WF {\text{WF}}

\def \beqq {\begin{equation}}
\def \eeqq {\end{equation}}
\def \bpf {\begin{proof}}
\def \epf {\end{proof}}
\def \beq {\begin{equation*}}
\def \eeq {\end{equation*}}

\def \La {\Lambda}    
   
\def \lap {\Delta}
\def \p {\partial}
\def \ha {\frac{1}{2}}   

\def \mck {{\mathcal K}}

\def \mcx {{\mathcal X}}

\def \mbn {{\mathbb N}}
\def \mbr {{\mathbb R}}
\def \mbs {{\mathbb S}}

\def \msf {{\mathscr F}}
\def \mss {{\mathscr S}}

\DeclareMathOperator*{\argmin}{arg\,min}                   
\renewcommand{\t} {^{\top}}                                

\renewcommand{\phi}{\mathbf{\varphi}}

\newcommand{\calK}{\mathcal{K}}

\newcommand{\bfA}{\mathbf{A}}

\newcommand{\bfb}{\mathbf{b}}

\newcommand{\bfX}{\mathbf{X}}
\newcommand{\bfx}{\mathbf{x}}

\newcommand{\bfe}{\mathbf{e}}

\newcommand{\bfy}{\mathbf{y}}

\newcommand{\bfM}{\mathbf{M}}

\newcommand{\bfV}{\mathbf{V}}

\newcommand{\bfQ}{\mathbf{Q}}

\newcommand{\bbR}{\mathbb{R}}

\def \mck {{\mathcal K}}

\def \mcx {{\mathcal X}}

\def \mbn {{\mathbb N}}
\def \mbr {{\mathbb R}}
\def \mbs {{\mathbb S}}

\def \msf {{\mathscr F}}
\def \mss {{\mathscr S}}
 
\def \id {\operatorname{Id}}
\def \comp {\operatorname{comp}}
\def \loc {\text{loc}}

\def \re {\operatorname{Re}}

\def \WF {\text{WF}}

\def \beqq {\begin{equation}}
\def \eeqq {\end{equation}}
\def \bpf {\begin{proof}}
\def \epf {\end{proof}}
\def \beq {\begin{equation*}}
\def \eeq {\end{equation*}}

\def \La {\Lambda}    
   
\def \lap {\Delta}
\def \p {\partial}
\def \ha {\frac{1}{2}}

\newcommand\ChangeRT[1]{\noalign{\hrule height #1}}

\begin{document}

\maketitle

\begin{abstract} 
We consider the imaging of cosmic strings by using Cosmic Microwave Background (CMB) data. Mathematically, we study the inversion of an X-ray transform in Lorentzian geometry, called the light ray transform. The inverse problem is highly ill-posed, with additional complexities of being large-scale and dynamic, with unknown parameters that represent multidimensional objects.  This presents significant computational challenges for the numerical reconstruction of images that have high spatial and temporal resolution.  In this paper, we begin with a microlocal stability analysis for inverting the light ray transform using the Landweber iteration. 
Next, we discretize the spatiotemporal object and light ray transform and consider iterative computational methods for solving the resulting inverse problem.  We provide a numerical investigation and comparison of some advanced iterative methods for regularization including Tikhonov and sparsity-promoting regularizers for various example scalar functions with conormal type singularities.
\end{abstract}



\section{Introduction} \label{sec-intro}

Our study is motivated by recent advances on the inverse problem of recovering gravitational anomalies such as cosmic strings from the Cosmic Microwave Background (CMB) radiation measurements. Since its discovery by Penzias and Wilson in 1965, the CMB has provided invaluable information regarding the early universe. Despite being highly isotropic, there are faint anisotropies in the CMB which have been mapped by sensitive detectors such as COBE, WMAP, and Planck Surveyor. It has been known for a long time \cite{sachs1967perturbations} that on the first-order linearization level, the CMB anisotropy is connected to  gravitational perturbations via a cosmological X-ray transform, also called the light ray transform. 
Recently, the possibility of recovering such gravitational perturbations from the integral transform received lots of attention; see \cite{LOSU, LOSU1, VaWa, Wan, Wan1}.  The mathematical studies are quite encouraging; however, research on the numerical inversion of the light ray transform remains preliminary. The primary goal of this work is to investigate the numerical inversion of the light ray transform, especially for recovering cosmic strings which are represented by singularities in the gravitational perturbations. We remark that the light ray transform naturally appears in other applications; see e.g.\ Section 5 of \cite{LOSU1}.

\begin{figure}[htbp]
\centering
\includegraphics[scale = 0.6]{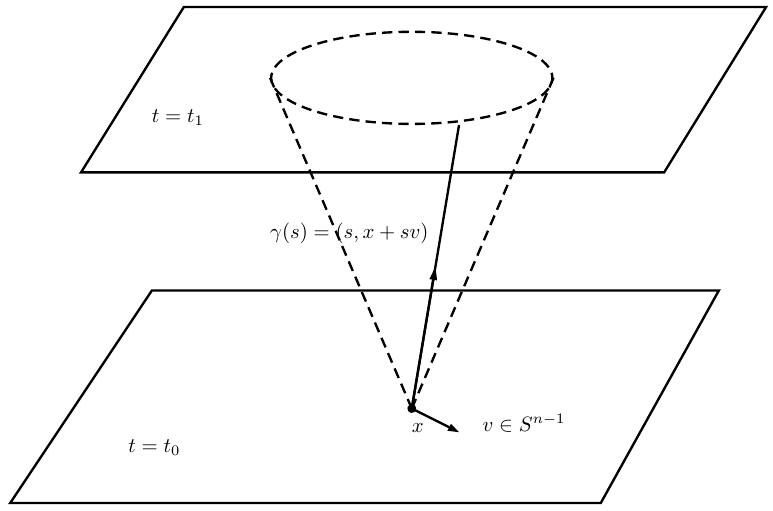}
\caption{Setup of the inverse problem. Each point at $t = t_0$ plane is regarded as a light source. The light signal from $x$ is recorded at $t = t_1$ in the dashed circle. The light ray $\gamma$ can be parametrized by $x\in \mbr^n, v\in \mbs^{n-1}$. The inverse problem is to reconstruct $f$ supported between $t = t_0$ and $t = t_1$ from $Lf$. }
\label{fig-setup}
\end{figure}

We begin with a brief description of the inverse problem. Consider the Minkowski spacetime $(\mbr^{n+1}, g), n\geq 2$ where in local coordinates $z = (t, x), t\in \mbr, x\in \mbr^n$,  $g = -dt^2 + \sum_{i = 1}^n (dx^i)^2$. We describe all null geodesics (or light-like geodesics) on $(\mbr^{n+1}, g)$ as
\beq
\gamma_{x, v}(s) = (s, x + sv), \quad s\in \mbr, x\in \mbr^n, v\in \mbs^{n-1}. 
\eeq
In particular, the set of null geodesics can be identified with $\mbr^n\times \mbs^{n-1}$.  The light ray transform for scalar functions can be defined as 
\beqq\label{eq-lray}
Lf(x, v) = \int_\mbr f(s, x+ sv)ds
\eeqq
where $f$ is a function so that the integral makes sense. The  inverse problem is to reconstruct $f$ from $Lf$. See Figure \ref{fig-setup}. 

The integral transform $L$ can be regarded as a Lorentzian analogue of the famous Radon transform which is the backbone of X-ray computed tomography; see for instance \cite{Nat}. Perhaps the more precise analogue is the limited angle Radon transform, because only integrations along null geodesics are considered. Due to the limited available data, the inversion of $L$ is severely ill-posed. As shown in several theoretical works \cite{Gui, LOSU, Wan}, not all information of $f$ can be stably recovered, even though the transform is injective in some cases, see for instance \ \cite{FIO, LOSU1, Ste}. Thus for numerical reconstruction, it is crucial to understand the mechanism of instability and find ways to improve it. Additionally, the discretized problem is highly ill-posed and regularization strategies have yet to be considered.

\paragraph{Overview of contributions}
In this work, we develop numerical methods for inverting the light ray transform in relation to imaging cosmic strings in a two-dimensional universe over time from the CMB measurements.  We investigate various iterative reconstruction methods and analyze the artefact phenomena. Our main contributions include the following.
\begin{itemize}
\item In a functional regime, we analyze stability and instability of the Landweber iteration method for inverting the light ray transform and use microlocal analysis to investigate the artefact phenomena. 
\item We discretize the linear inverse problem, formulating it as a large, sparse matrix equation and consider various regularization techniques for the reconstruction of cosmic strings in $\mathbb{R}^{2+1}$. 
\item We provide a numerical comparison of various iterative methods, including iterative reconstruction methods like Landweber, iterative Krylov projection methods for Tikhonov regularization, and iterative shrinkage threshholding algorithms for $\ell_1$ regularization.  The computational methods that we use are state-of-the-art, but they are not new. For the tomography approach for the CMB inverse problems, these methods are only now becoming relevant.  
\end{itemize}

An outline of the paper is as follows.  In Sections \ref{sec-land} and \ref{sec-art}, we focus on the functional representation of the linearized CMB inverse problem. We describe the Landweber iteration and analyze its convergence behavior in an unstable regime.  Then in Section \ref{sec-comp} we describe the discrete inverse problem and describe various iterative methods that can be used for the reconstruction of cosmic strings. Several numerical examples are provided in Section \ref{sec-results} to illustrate the theory and demonstrate the potential benefits and challenges of incorporating different regularization terms.  Conclusions and future work are provided in Section \ref{sec-conc}.
\section{The Landweber Iteration}\label{sec-land}
We recall  from Section 3 of \cite{StYa} the Landweber's iteration in the abstract setting. Let $A: H_1 \rightarrow H_2$ be a linear operator where $H_1, H_2$ are Hilbert spaces. For a given $m\in H_2$, we aim to solve the linear inverse problem $Af = m$ for $f\in H_1$. It is well-known that Landweber iteration is a gradient descent method for the functional 
\beq
f \rightarrow \|Af - m\|_{H_2}^2. 
\eeq
We can derive it from solving $Af = m$  in the form 
\beq
(\id - (\id - \gamma A^*A))f = \gamma A^*m
\eeq
where $A^*$ is the adjoint of $A$, $\id$ is the identity operator and $\gamma>0$ is chosen properly. Let $K = \id - \gamma A^*A$. If $K$ is a strict contraction, then $f$ can be solved by Neumann series, 
\beq
f = \sum_{j = 0}^\infty K^j \gamma A^*m.
\eeq
This implies the following Landweber iteration scheme: 
\beqq\label{eq-land}
\begin{gathered}
f_0 = 0 \text{ and } f_k = f_{k - 1} - \gamma A^*(A f_{k-1} - m), \quad k = 1, 2, \ldots. 
\end{gathered}
\eeqq
Convergence of the iteration in the abstract setting was analyzed in Section 3 of \cite{StYa}. For example, if the stability estimate
\beqq\label{eq-stab0}
\mu \|f\|_{H_1}\leq \|Af\|_{H_2}
\eeqq
holds for some $\mu>0$, then $K$ is a strict contraction for $\gamma \in (0, 2/\|A\|^2)$. In fact, the Neumann series converges uniformly. 
Even when $A$ has a non-trivial kernel,  convergence can be established under a stability assumption similar to \eqref{eq-stab0}, see Theorem 3.1 of \cite{StYa}.  

When applying the Landweber iteration to the light ray transform $L$, a serious issue is that \eqref{eq-stab0} does not hold. For simplicity, let $\mck$ be a pre-compact set of $\mbr^{n+1}$ and we consider $f$ supported in $\mck$. Then there is no stability estimate of the form 
 \beq
 \|f\|_{H^a}\leq C\|Lf\|_{H^{b}}, 
 \eeq
 for any $a, b\in \mbr$, $C>0$. Similar behavior is known for the limited angle Radon transform, see Section VI.2 of \cite{Nat}. We remark that for $f\in L^2$ with compact support, $Lf = 0$ implies $f = 0$, see the discussion below. So the kernel of $L$ is trivial and the instability is not caused by the kernel.

To further understand the issue, we discuss the injectivity and stability of the light ray transform. We will obtain a partial stability estimate and use it to show the convergence of the Landweber iteration. 
We start from the injectivity proof of $L$, see for example \cite{LOSU1}.  Suppose $f\in C_0^\infty(\mbr^{n+1})$ and $Lf(x, v) = 0$ for all $x\in \mbr^n, v\in \mbs^{n-1}.$ We denote the Fourier transform on $\mbr^{n+1}$ by $\msf$ and the Fourier transform in the spatial variable by $\overline\msf$. Taking the Fourier transform of \eqref{eq-lray} in the $x$ variable, we get 
\beqq\label{eq-fourier}
\begin{aligned}
\overline\msf(L f)(\xi, v) &= \int_{\mbr^n} \int_\mbr e^{-\imath x\cdot \xi} f(t, x+tv) dt dx \\
& = \int_{\mbr^n} \int_\mbr e^{-\imath y \cdot \xi - \imath t(v \cdot \xi)} f(t, y) dt dy = \msf f(v\cdot \xi, \xi)
\end{aligned}
\eeqq
where $\imath^2 = -1.$ 
This is known as the {\em Fourier Slice Theorem}. Note that for $v\in \mbs^{n-1}, \xi\in \mbr^n$, $(v\cdot\xi, \xi)$ gives all space-like and light-like vectors on $(\mbr^{n+1}, g)$. Thus if $Lf = 0$, we get $\msf f= 0$ in the cone of  space-like vectors. By analyticity of $\msf f$, we conclude that $\msf f = 0$ on $\mbr^{n+1}$ hence $f = 0.$ Indeed, we see that the kernel of $L$ acting on Schwartz functions $\mss(\mbr^{n+1})$ consists of $f$ whose Fourier transform vanishes on all space-like vectors. 

The above argument is clearly not stable. However, there is some stable ingredient in the argument. We denote by  $\Gamma^{sp} = \{(z, \zeta)\in T^*\mbr^{n+1}: \zeta = (\tau, \xi), \tau \in \mbr, \xi\in \mbr^n, |\xi| > |\tau|\}$ the cone of space-like (co)vectors. For $\delta \geq 0$ small, we let $\Gamma^{sp}_\delta = \{(z, \zeta)\in T^*\mbr^{n+1}: \zeta = (\tau, \xi), \tau \in \mbr, \xi\in \mbr^n, (1 - \delta)|\xi| > |\tau|\}$ so $\Gamma^{sp}_\delta\subset \Gamma^{sp}$. Let $\chi_\delta(\tau, \xi)$ be the characteristic function of $\Gamma^{sp}_\delta$ and let $\chi_\delta(D) = \msf^{-1}\chi_\delta\msf$ be defined as a Fourier multiplier.  Below, we use the product measure to define the $L^2$ space on $\mbr^n\times \mbs^{n-1}.$ 
\begin{prop}\label{prop-stab}
For any $f\in L^2(\mbr^{n+1}), n\geq 2$, we have
\beqq\label{eq-stab1}
\mu_n \|\chi_\delta(D) f\|_{L^2}\leq  \|Lf\|_{H^{1/2}}.
\eeqq
where $\mu_n = C_n^\ha \delta^{(n-3)/4}>0$ and $C_n = 2\pi^{(n-1)/2}/\Gamma((n-1)/2)$ is the surface area of $\mbs^{n-2}, n\geq 2$.  
\end{prop}
\bpf
From \eqref{eq-fourier}, we get 
\beq
\begin{gathered}
\overline\msf(L f)(\xi, v) =   \chi_0(v\cdot\xi, \xi) \msf f(v\cdot \xi, \xi)
\end{gathered}
\eeq
almost everywhere. By using the Plancherel theorem, we obtain that 
\beq
\begin{gathered}
\|L f\|^2_{H^{1/2}} \geq \int_{\mbr^n}\int_{\mbs^{n-1}}|\xi| |\chi_0(v\cdot\xi, \xi) \msf f(v\cdot \xi, \xi) |^2 dv d\xi.
\end{gathered}
\eeq
To compute the integral in $v$, without loss of generality, we can take $\xi = |\xi|(1, 0, 0)$ and write $v = (\tau, \omega)$ where $\tau\in [-1, 1]$ and $\omega(\sqrt{1-|\tau|^2})^{-1}\in \mbs^{n-2}$. Then we get 
\beq
\begin{aligned}
\|L f\|^2_{H^{1/2}}  &\geq \int_{\mbr^n}\int_{\mbs^{n-2}} \int_{-1}^1|\xi|  |\chi_0(\tau |\xi|, \xi) \msf f(\tau |\xi|, \xi)|^2 d\tau d\omega d\xi  \\
&= C_n \int_{\mbr^n}  \int_{-1}^1 (\sqrt{1-|\tau|^2})^{n-3}|\xi|  |\chi_0(\tau |\xi|, \xi) \msf f(\tau |\xi|, \xi) |^2   d\tau   d\xi \\
&=  C_n \int_{\mbr^n}  \int_{-|\xi|}^{|\xi|} \frac{(|\xi|^2-|\tau|^2)^{(n-3)/2}}{|\xi|^{n-3}}  |\chi_0(\tau , \xi) \msf f(\tau , \xi) |^2   d\tau   d\xi 
\end{aligned}
\eeq 
Using the support property of $\chi_0,$ we have
\beqq\label{eq-l0f}
\begin{aligned}
\|L f\|^2_{H^{1/2}}   
 &\geq  C_n  \int_{\mbr^n}  \int_{-|\xi|}^{|\xi|} \frac{(|\xi|^2-|\tau|^2)^{(n-3)/2}}{|\xi|^{n-3}} |\chi_\delta(\tau , \xi) \msf f(\tau , \xi) |^2   d\tau   d\xi \\
  &\geq  C_n  \int_{\mbr^n}  \int_{-|\xi|}^{|\xi|} \frac{(\delta|\xi|^2)^{(n-3)/2}}{|\xi|^{n-3}} |\chi_\delta(\tau , \xi) \msf f(\tau , \xi) |^2   d\tau   d\xi \\
 & \geq C_n\delta^{(n-3)/2}  \|\chi_\delta(D) f\|^2_{L^2}
\end{aligned}
\eeqq  
This completes the proof. 
 \epf
 
Based on the above results, we can set up the Landweber iteration to stably recover certain information of $f$. We aim to solve $Lf = m$ for $f\in L^2$ supported in $\mck$. By the description of the kernel of $L$, we see that the equation is equivalent to $L\chi_0(D)f = m.$ By the continuity of $L: L^2(\mbr^{n+1})\rightarrow H^{1/2}(\mbr^n\times \mbs^{n-1})$ (see for example \cite{Wan}), we know that the equation should be solved in $H^{1/2}$. As in \cite{HMS}, it is convenient to further modify the equation to $L^2.$ Let $L^*$ be the $L^2$ adjoint. We can solve $L^* L \chi_0(D)f = L^* m$ in $H^{1}$. To get a stable reconstruction, we consider solving $\chi_\delta(D) L^* L \chi_0(D)f = \chi_\delta(D) L^* m$ in $H^{1}$. Let $\lap$ be the Dirichlet realization of the positive Laplacian on $\mck$. We can define the square root of $\lap$. Then we solve 
\beqq\label{eq-inv}
\lap^\ha \chi_\delta(D) L^* L \chi_0(D)f = \lap^\ha \chi_\delta(D) L^* m
\eeqq 
in $L^2$ by Landweber iteration.  

Let $X = \lap^\ha \chi_\delta(D) L^* L\chi_\delta(D)$.  The corresponding Landweber iteration for solving \eqref{eq-inv} is 
\beqq\label{eq-land1}
\begin{gathered}
f_0 = 0 \text{ and } f_k = f_{k - 1} - \gamma X^*(X f_{k-1} - \lap^\ha L^*m), \quad k = 1, 2, \cdots. 
\end{gathered}
\eeqq
Note that $X: L^2(\mbr^{n+1})\rightarrow L^2(\mbr^{n+1})$ is bounded. Also, from Proposition \ref{prop-stab}, we have 
\beq
\mu_n^2\|\chi_\delta(D)f\|_{L^2}^2 \leq \|L\chi_\delta(D)f\|_{H^{1/2}}^2 \leq \|Xf\|_{L^2} \|\chi_\delta(D) f\|_{L^2}.
\eeq
Then $K= \id - \gamma  X^*X$ is a strict contraction for $\gamma\in (0, 2/\|X\|^2)$, and we conclude that
\beq
\chi_\delta(D)f = \sum_{j = 0}^\infty K^j \gamma \lap^\ha \chi_\delta(D) L^*m
\eeq
 converges uniformly. This means that $\chi_\delta(D) f$ can be stably reconstructed from $Lf.$

\section{Analysis of the Instability}\label{sec-art}
In this section, we analyze the behavior of the Landweber reconstruction in the unstable regime. We are particularly interested in the artefact phenomena. Note that the iteration \eqref{eq-land1} is basically for solving $Lf = m$ in $H^{1/2}$. This is suitable for proving convergence of the scheme. However, as we will see soon, the iteration will produce stronger artefacts and it might be preferable to solve $Lf = m$ in $L^2$ to reduce the artefacts.  We thus consider
\beqq\label{eq-land2}
\begin{gathered}
f_0 = 0 \text{ and } f_k = f_{k - 1} - \gamma L^*(L f_{k-1} - L f), \quad k = 1, 2, \ldots. 
\end{gathered}
\eeqq

Let $N = L^* L$ be the normal operator of $L$. It suffices to analyze the microlocal structure of $N$ and powers of $N$. For example, the first few terms of the iteration \eqref{eq-land2} are 
\beqq\label{eq-land3}
\begin{gathered}
f_0 = 0,  \quad f_1 =   \gamma  N  f, \quad  f_2 = 2\gamma N f  - \gamma^2 N^2  f.
\end{gathered}
\eeqq
 It is computed in \cite[Theorem 2.1]{LOSU1} that 
\beqq\label{eq-minnormal}
N f(t, x)  = \int_{\mbr^{n+1}} K_N(t, x, t', x') f(t', x') dt'dx'
\eeqq
where the Schwartz kernel is given by 
\beqq\label{eq-norker}
K_N(t, x, t', x') = \frac{\delta(t - t' - |x - x'|) + \delta(t - t' + |x - x'|)}{|x - x'|^{n - 1}}.
\eeqq 
In particular, $N$ can be written as a Fourier multiplier 
\beqq\label{eq-normalfourier}
\begin{gathered}
N f(t, x) =   \int_{\mbr^{n+1}} e^{i(t \tau +x\cdot \xi)} k(\tau, \xi) \hat f(\tau, \xi) d\tau d\xi
\end{gathered}
\eeqq
where 
\beqq\label{eq-ksym}
k(\tau, \xi) = C_n \frac{(|\xi|^2 - \tau^2)_+^{\frac{n - 3}{2}}}{|\xi|^{n-2}}, \quad C_n = 2\pi |\mbs^{n-2}|.
\eeqq
Here, for $s\in \mbr$ and $\re a> -1$, $s_+^a$ denotes the  distribution defined by $s_+^a = s^a$ if $s > 0$ and $s_+^a = 0$ if $s \leq 0$, see \cite[Section 3.2]{Ho1}. Furthermore, 
\beqq\label{eq-Nf}
N  f(t, x) =   \int_{\mbr^{n+1}} \int_{\mbr^{n+1}} e^{i(t - t', x - x')\cdot (\tau, \xi)} k(\tau, \xi)  f(t', x') d\tau d\xi dt'dx'. 
\eeqq
Note that the amplitude $k$ has conormal type singularities at $\tau^2 = |\xi|^2$. So the operator is not a standard pseudo-differential operator. When restricting $k(\tau, \xi)$ to $\Gamma^{sp}$, we indeed get an elliptic pseudo-differential operator which is the main theme in \cite{LOSU, LOSU1}. In particular, this means that one can recover space-like singularities in $f$ from $Nf$. For time-like directions, $k(\tau, \xi) = 0$ thus time-like singularities of $f$ are lost in $Nf$.  For the borderline case of light-like singularities, we need the microlocal structure of $N$ in \cite{Wan} that the Schwartz kernel is a paired Lagrangian distribution. Such distributions were introduced in \cite{GuUh, MeUh}. We briefly recall the relevant result for understanding $N.$

Let $\mcx$ be a $C^\infty$ manifold of dimension $n$ and $w_\mcx$ be the simplectic form on $T^*\mcx$.  Let $\La_0, \La_1$ be conic Lagrangian submanifolds of $T^*(\mcx\times \mcx)\backslash 0$ with symplectic form $\pi_1^*w_\mcx + \pi_2^*w_\mcx$. Here, $\pi_1, \pi_2: \mcx\times \mcx\rightarrow \mcx$ denotes the projections to the first, second copy of $\mcx$.  Suppose that $\La_1$ intersects $\La_0$ cleanly at a codimension $k$, $1\leq k\leq 2n-1$ submanifold $\Sigma = \La_0\cap \La_1$, namely
\beq
T_p(\La_0\cap \La_1) = T_p(\La_0)\cap T_p(\La), \quad \forall p\in \Sigma. 
\eeq
 
It is proved in Theorem 3.1 of \cite{Wan} that the Schwartz kernel of the normal operator $K_N \in I^{-n/2, n/2- 1}(\mbr^{n+1}\times \mbr^{n+1}; \La_0, \La_1)$, in which $\La_0, \La_1$ are two cleanly intersecting Lagrangians  defined as follows: 
 \beqq\label{eq-lag1}
 \begin{gathered}
 \La_0 = \{(t, x, \tau, \xi; t', x', \tau', \xi')\in T^*\mbr^{n+1}\backslash 0 \times T^*\mbr^{n+1}\backslash 0: \\
 t'  = t, x' = x, \tau' = -\tau, \xi' = -\xi\}
 \end{gathered}
 \eeqq
which is the punctured conormal bundle of the diagonal in $\mbr^{n+1}\times \mbr^{n+1}$ and 
 \beqq\label{eq-lag2}
 \begin{gathered}
 \La_1 = \{(t, x, \tau, \xi; t', x', \tau', \xi')\in T^*\mbr^{n+1}\backslash 0 \times T^*\mbr^{n+1}\backslash 0: x = x' + (t - t')\xi/|\xi|,\\
 \tau = \pm|\xi|,  \tau' = -\tau,  \xi' = -\xi \}.
  \end{gathered}
 \eeqq
 
To analyze the powers of $N$, we use the fact that $(\La_0, \La_1)$ is the flow-out model, see \cite{Wan}. We recall that a submanifold $\Gamma \subseteq T^*\mcx$ is involutive if $\Gamma = \{(x, \xi): p_i(x, \xi) = 0, i = 1, 2, \cdots, k\}$ satisfies (i) $p_i$ are defining functions of $\Gamma$  that $dp\neq 0$ on $p = 0$,  and (ii) $p_i$ are in involution so the Poisson brackets $\{p_i, p_j\} = 0$ at $\Gamma$. Let $H_{p_i}$ be the Hamilton vector fields of $p_i$. The flow out of $\Gamma$ 
\beqq\label{eq-flowout}
\begin{gathered}
\La_\Gamma = \{(x, \xi, y, \eta)\in T^*\mcx \times T^*\mcx: (x, \xi)\in \Gamma,  \\
(y, \eta) = \exp(\sum_{j = 1}^k t_j H_{p_j})(x, \xi), t_j \in \mbr \}
\end{gathered}
\eeqq
is a Lagranian submanifold of $T^*(\mcx \times \mcx)$ and is a canonical relation if $\Gamma$ is conic. Write $\La_1 = \La_\Gamma$. In this case, we have the following result.
\begin{theorem}[Theorem 1 of \cite{AnUh}]\label{thm-calculus}
Let $A\in I^{p, l}(\mcx\times \mcx; \La_0, \La_1)$ and  $B\in I^{r, s}(\mcx\times \mcx; \La_0, \La_1)$ be properly supported. Then $A\circ B\in  I^{p+ r + k/2, l + s - k/2}(\mcx\times \mcx; \La_0, \La_1)$. 

\end{theorem}

Now consider the two Lagrangians $\La_0$ and $\La_1$ in \eqref{eq-lag1} and \eqref{eq-lag2} respectively. Let $f(\tau, \xi) = \ha(\tau^2 - |\xi|^2)$ and $\Sigma = \{(t, x, \tau, \xi; t, x, -\tau, -\xi)\in T^*\mbr^{n+1}\backslash 0 \times T^*\mbr^{n+1}\backslash 0: f(\tau, \xi) = 0\}$.  Then $\La_1$ is the flow out of $\Sigma$ under the Hamilton vector field $H_f$. Indeed, we have 
\beq
H_f = \tau \frac{\p}{\p t} - \sum_{i = 1}^3 \xi_i \frac{\p}{\p x^i}.
\eeq
Let $\gamma(s) = (t(s), x(s), \tau(s), \xi(s))$ be a null bi-characteristic which satisfies
\beq
\begin{gathered}
\dot t(s) = \tau, \quad \dot x_i(s) = -\xi_i, \quad \dot \tau(s) = 0, \quad \dot \xi_i(s) = 0\quad s\geq 0\\
 t(0) = t', \quad x(0) = x', \quad \tau(0) = \tau', \quad \xi(0) = \xi'
\end{gathered}
\eeq
with $f(\tau', \xi') = 0.$ We solve 
\beq
t(s) = t' +  s \tau', \quad x(s) = x' - s \xi', \quad \tau(s) = \tau', \quad \xi(s) = \xi',  
\eeq
which up to a re-parametrization gives \eqref{eq-lag2}. Now we can use Theorem \ref{thm-calculus} to conclude that for $j\in \mbn$, $N^j \in  I^{-nj/2 + (j-1)/2, nj/2 - (j-1)/2}(\mbr^{n+1}\times \mbr^{n+1}; \La_0, \La_1)$. So all the terms $N^j, j = 1, 2, \ldots$ have similar microlocal structures but the orders are different. 
 
It is known that for a paired Lagrangian distribution $N^j$, the wave front set $\WF(N^j)\subset \La_0\cup \La_1$. From wave front analysis see e.g.\ \cite{Ho1}, we observe that the wave front set on $\La_0$ does not produce new singularities in $N^jf$. However, the part on $\La_1$ does produce new singularities which correspond to the artefacts in the numerical reconstruction. In fact, away from $\La_0$, we know that $N^j \in I^{-(n-1)j/2 - 1/2}(\mbr^{n+1}\times \mbr^{n+1}; \La_1)$ is a Lagrangian distribution. By the definition of $\La_1$ and wave front analysis, we see that only light-like singularities are affected by $\La_1$, and the new singularities (or artefacts) would appear in the flow-out  under $\La_1$. The analysis so far implies that in the iteration scheme \eqref{eq-land2}, one can reconstruct space-like singularities of $f$ but not time-like singularities, and there might be artefacts from light-like singularities of $f$.

 Below, we give a more precise statement about the strength of the artefacts. As in \cite{Qui}, we describe them using the notion of microlocal Sobolev regularity.  We recall from Definition 2.2 of \cite{Qui} that a distribution $f$ is in $H^s$ microlocally near $(z_0, \zeta_0)\in \mbr^{n+1}\times \mbr^{n+1}\backslash 0$ if and only if there exists a cut-off function $\phi\in C_0^\infty(\mbr^{n+1})$ with $\phi(z_0)\neq 0$ and $u(\zeta)$ homogeneous of degree zero and smooth on $\mbr^{n+1}\backslash 0$ with $u(\zeta_0)\neq 0$ such that $(1 + |\zeta|^2)^{s/2}u(\zeta)\msf (\phi f)(\zeta)\in L^2(\mbr^{n+1})$. It follows from the definition that if $(z_0, \zeta_0)\notin \WF(f)$, then $f$ is $H^s$ microlocally near $(z_0, \zeta_0)$ for all $s\in \mbr.$
 
 \begin{theorem}\label{thm-sobo}
Suppose $f$ is $H^s$ microlocally near $(z, \zeta) = (t, x, \tau, \xi)$ with $\tau^2 = |\xi|^2$. 
Then for $j\in \mbn$, $N^j f$ is $H^{s + (n-1)j/2}$ microlocally near $(z', \zeta')$ where $z' = (t', x - (t' - t)\xi/|\xi|)$ and $\zeta = (\tau, \xi)$.
 \end{theorem}
 \bpf
 Let $C_1 = \La_1'$ be the canonical relation. We can verify that the differential of the projection $C_1 \rightarrow T^*\mbr^{n+1}$ has rank $2n + 1$. Thus applying Theorem 4.3.2 of \cite{Ho0}, we get that $N^j: L^2_{\comp}(\mbr^{n+1})\rightarrow H^{(n-1)j/2}_{\loc}(\mbr^{n+1})$.  Then we can follow the arguments in Theorem 3.1 of \cite{Qui} to get the microlocal result. 
 \epf
 
We note that as $j$ increases, the produced artefacts are more regular in Sobolev regularity. For the iteration scheme \eqref{eq-land1}, the corresponding artefacts are actually more singular after each iteration due to the $\lap^\ha$ terms. This is why we prefer the iterative scheme \eqref{eq-land2}. 

\begin{figure}[htbp]
\centering
\includegraphics[scale = 0.6]{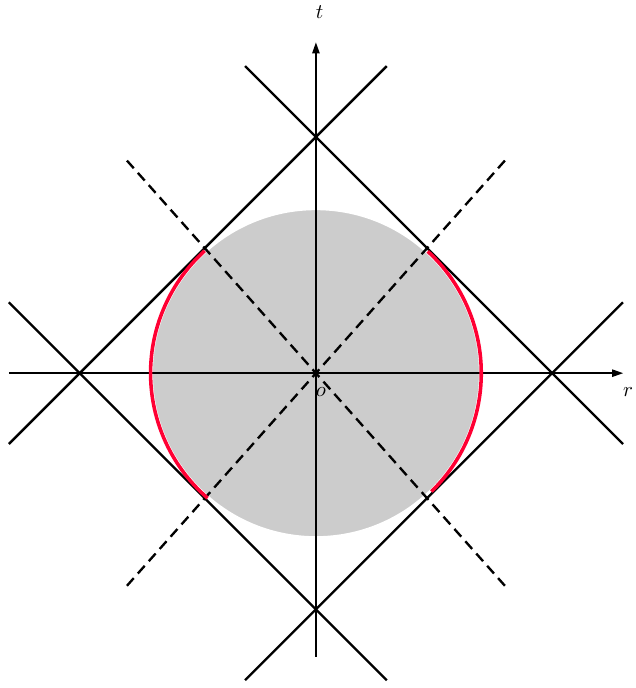}
\caption{Reconstruction of the characteristic function of the unit ball. The figure shows the projection to the $t, x$-plane so that each point represents a circle in $\mbr^2.$  On the circle in the figure, the highlighted arches (in red) represent parts of the sphere that can be stably reconstructed. The four straight lines tangent to the sphere represent possible artefacts in the reconstruction.  }
\label{fig-sphere}
\end{figure}

Finally, we use one example to demonstrate all the possible phenomena in the reconstruction based on the analysis.  Consider $n=2$ and let $f$ be the characteristic function of the unit ball in $\mbr^3$, that is for $t, x, y\in \mbr$, 
\beq
f(t, x, y) = 1, \quad \text{ if } t^2 + x^2 + y^2 < 1
\eeq 
and otherwise $f = 0.$ In Figure \ref{fig-sphere}, we show the projection of $f$ to the $t, x$-plane. The function has conormal type singularities at the unit sphere. We can identify the part of the boundary that can be reconstructed stably and the possible artefacts. 

\section{Computational Approaches for the Inverse Problem} 
\label{sec-comp}
We turn our attention to computational approaches for solving the linear inverse problem in a discrete setting. We begin in Section \ref{sec-forwardProb} with the forward problem, describing how the space-time is discretized and how the forward light ray transform model can be represented using a large, sparse matrix. Section \ref{sec-iterMeth} then provides an overview of iterative methods that can be used for string recovery.  We consider the discrete Landweber method as an iterative regularization strategy, as well as variational regularization approaches, namely Tikhonov and $\ell_1$ regularization. Numerical investigations are provided in Section \ref{sec-results}. 

\subsection{Discrete Forward Problem} \label{sec-forwardProb}
The discrete inverse problem can be stated as follows. We consider a $2+1$ spatial-temporal region, where space has been discretized using equally spaced $xy$-spatial grids and time has been discretized using $T$ time points.  We denote $\textbf{X}(t) \in \mathbb{R}^{n \times n}$ to be the unknown (discretized) image plane at time point $t$ where $t=1,2,\dots,T$. Upon vectorizing each spatial slice columnwise so that
\begin{equation*}
    \text{vec}\Big(\mathbf{X}(t)\Big) = \mathbf{x}(t) \in \mathbb{R}^{n^2},
\end{equation*}
the vector of unknowns is given by 
\begin{equation*}
    \bfx = \begin{bmatrix} \bfx(1)  \\ \bfx(2) \\ \vdots \\ \bfx(T)\end{bmatrix} \in \bbR^{n^2T}.
\end{equation*}
Assume there are $N$ source grid locations from which a light-cone emanates. For each source location $s_i, i = 1, 2, \ldots, N$, there are $m_i$ observations that are collected in vector
\begin{equation} \label{localRayTrace}
\bfb_i = \bfA_i \bfx + \mathbf{e}_i
\end{equation}
where $\bfb_i \in \bbR^{m_i}$ and $\bfA_i \in \bbR^{m_i \times n^2T}$ with noise contamination $\bfe_i \in \mathbb{R}^{m_i}$. The values $m_i$ correspond to the number of observable detector grid locations such that light-rays from source point $s_i$ have a non-zero interaction with the spatial-temporal object.

\begin{figure}
 \centering
        \includegraphics[width = .5\textwidth]{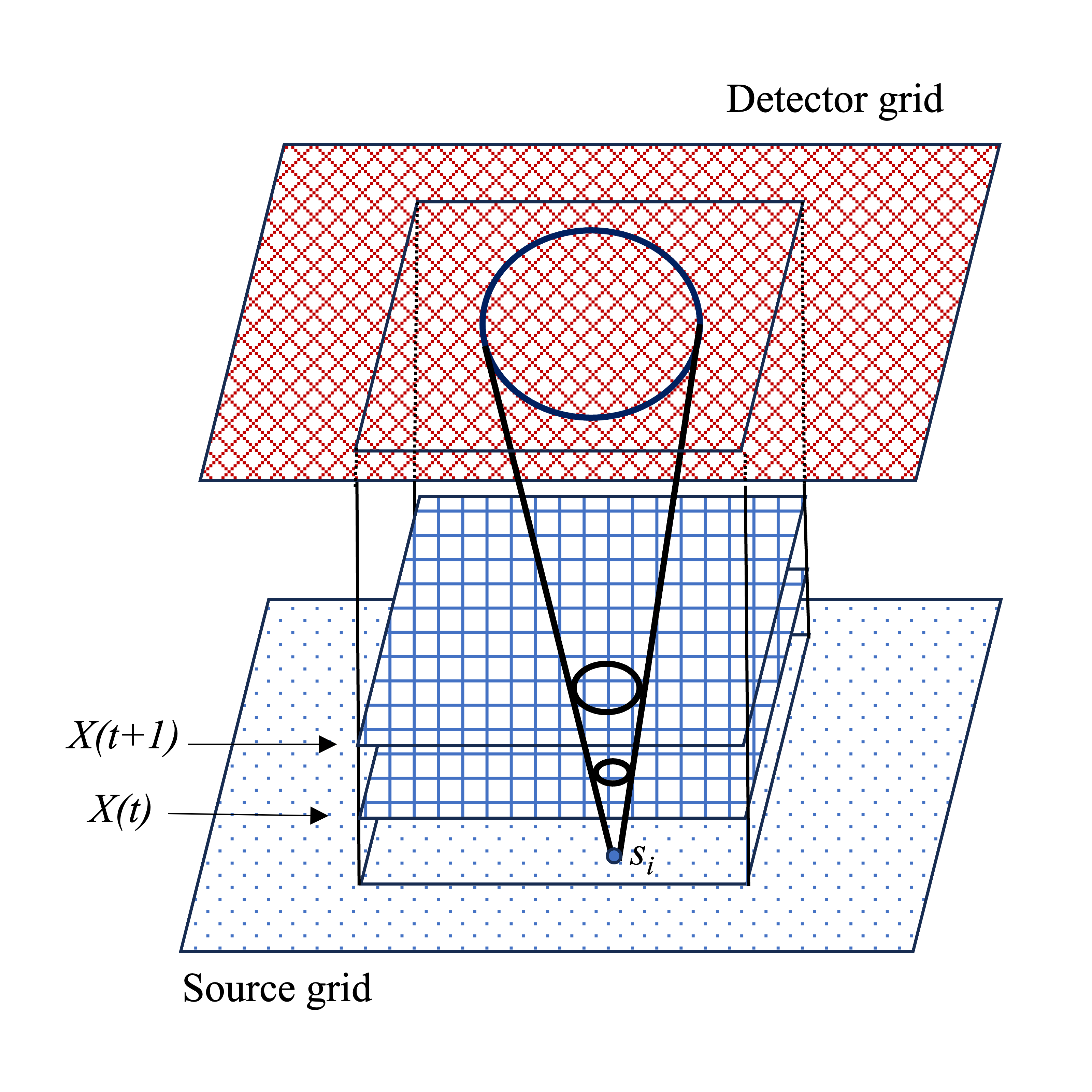}
        \includegraphics[width = .45\textwidth]{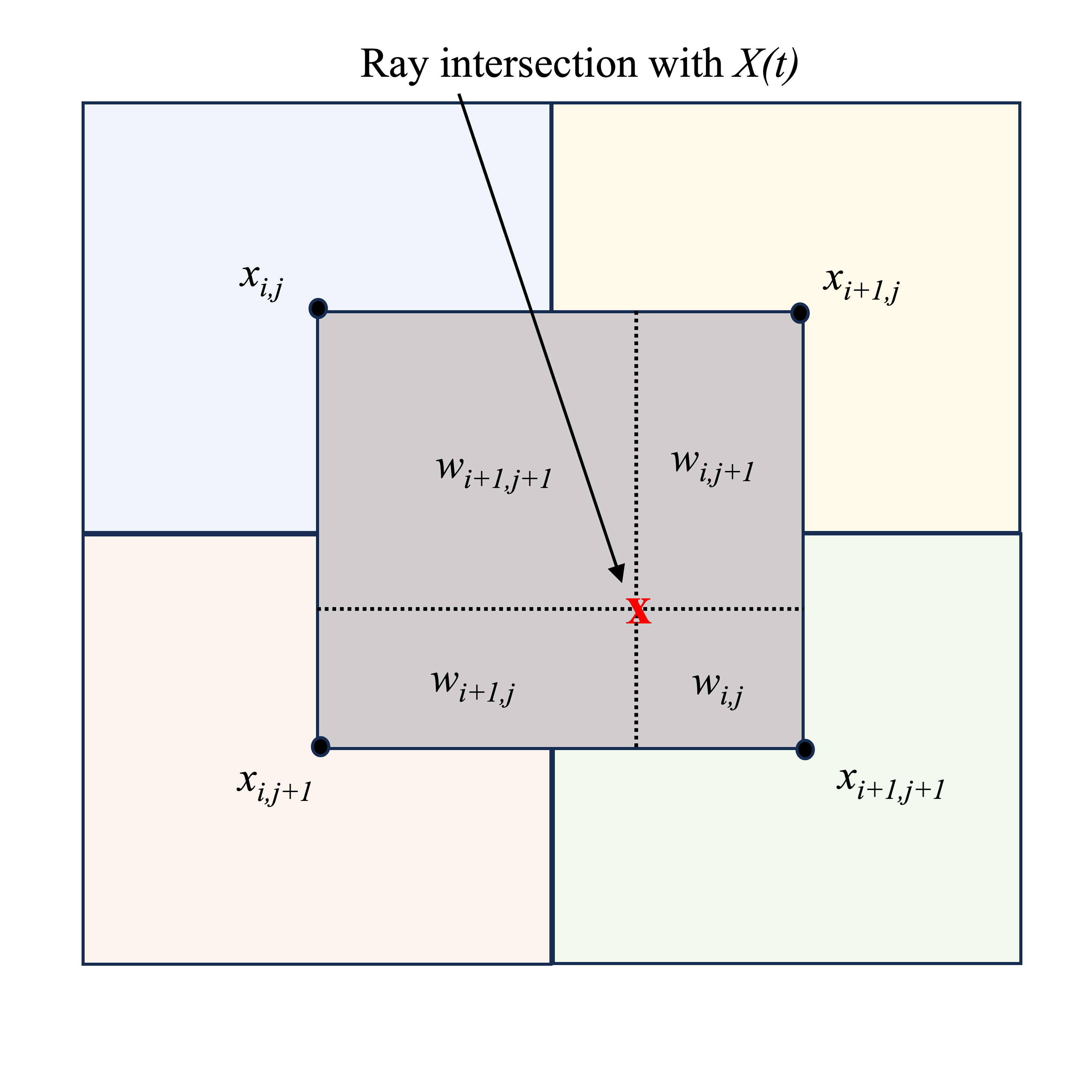}
\caption{Discrete problem setup.  On the left is an illustration of the geometry of the discrete CMB problem, with the object of interest located between the source grid at the bottom and the detector grid on the top.  For each source location $s_i$, detectors (determined by a cone of 45 degrees) receive contributions from event planes $X(t), t=1, \ldots, T$.  The weights of contribution from $X(t)$ are determined via bilinear interpolation, as shown in the right figure, where the red cross denotes the ray's intersection with $X(t)$. }
\label{discreteproblem}
\end{figure}

Equation \eqref{localRayTrace} models the linearized CMB problem outlined in Section \ref{sec-intro}. Here, the rows of $\mathbf{A}_i$ contain the weights of the intersection of the light-rays through the unknown spatial slices $\mathbf{X}(t)$ corresponding to source point $s_i$. The weights are obtained through bilinear interpolation.  More specifically, for each ray connecting the source $s_i$ with one detector, we find the point of intersection of the ray with the plane $\mathbf{X}(t)$ and identify the four nearest pixels $x_{i,j},\, x_{i+1,j},\, x_{i,j+1},$ and $x_{i+1,j+1}$.  Then the bilinear interpolation weights $w_{i,j},\, w_{i+1,j},\, w_{i,j+1},$ and $w_{i+1,j+1}$ are computed and placed in the corresponding columns of matrix $\bfA_i(t)$. See Figure \ref{discreteproblem} for an illustration.
Since the observations may contain contributions from multiple event planes, let $\bfA_i(t) \in \bbR^{m_i \times n^2}$ denote the ray-trace matrix that represents contributions from the plane $\bfX(t)$ when the source is located at $s_i$. Then we have for source point $s_i,$ 
\begin{equation*}
\bfA_i = \begin{bmatrix} \bfA_i(1) & \bfA_i(2) & \cdots & \bfA_i(T)\end{bmatrix} \in \bbR^{m_i \times n^2T}.
\end{equation*}

If we consider all $N$ source grid locations, we obtain the following linear inverse problem,
\begin{align} \label{DiscreteLinearProb}
\underbrace{\begin{bmatrix}\bfb_1 \\ \bfb_2 \\ \vdots \\ \bfb_N \end{bmatrix}}_{\bfb}&  =  \underbrace{\begin{bmatrix}\bfA_1 (1) & \bfA_1 (2) & \cdots &\bfA_1 (T) \\ \bfA_2 (1) & \bfA_2 (2) & \cdots &\bfA_2 (T)\\ \vdots & \vdots & \ddots & \vdots\\ \bfA_N(1) & \bfA_N (2) & \cdots &\bfA_N (T) \end{bmatrix}}_{\bfA} \underbrace{\begin{bmatrix} \bfx(1)  \\ \bfx(2) \\ \vdots \\ \bfx(T)\end{bmatrix}}_{\bfx} + \underbrace{\begin{bmatrix}\bfe_1 \\ \bfe_2 \\ \vdots \\ \bfe_N \end{bmatrix}}_{\bfe}
\end{align}
where $\bfb\in \bbR^m$ and $\bfA \in \bbR^{m \times n^2T}$
with $m = \sum_{i=1}^N m_i$. As every $\mathbf{A}_i(t)\in \mathbb{R}^{m_i \times n^2}$ has at most $4$ column entries per row, we are able to construct and store $\bfA$ as a sparse matrix; this enhances the efficiency of computing matrix-vector multiplications with $\bfA$ and $\bfA\t$. 
Next we describe iterative methods to estimate $\mathbf{x}$, given $\bfb$ and $\bfA$.

\subsection{Iterative Methods and Regularization} \label{sec-iterMeth}
Consider an inverse problem where noisy data is given by $\bfb = \bfA \bfx_{\rm true} + \bfe$ for some $\bfx_{\rm true}\in \bbR^{n^2T}.$  We consider the least-squares (LS) problem, 
\begin{equation} \label{LSprob}
\min_\bfx \| \bfA \bfx - \bfb \|^2,
\end{equation}
where unless specified otherwise, $\|\cdot\|$ will denote the Euclidean norm. Since $\bf{b}$ is contaminated by error and $\bfA$ is ill-conditioned, the solution of \eqref{LSprob} is usually a poor approximation of $\mathbf{x}_{\rm true}$ due to the propagation of error into the computed solution.
However, a regularized solution can be obtained by applying an iterative method to \eqref{LSprob} and terminating the iterations before the propagation of noise occurs.

One such method is the discrete Landweber method whose iterates are given by
\begin{equation} \label{landweb}
    \mathbf{x}^{(k+1)} = \mathbf{x}^{(k)} + \omega\,\bfA^T\mathbf{r}^{(k)},
\end{equation}
where $\omega \in (0,2/\sigma_1^2)$ is a relaxation term with $\sigma_1$ the largest singular value of $\bfA$. In the absence of noise and in exact arithmetic, we expect the errors to be nonincreasing, i.e., $\|\mathbf{x}_{\rm true} - \mathbf{x}^{(k+1)}\| \leq \|\mathbf{x}_{\rm true} - \mathbf{x}^{(k)}\|$. However, due to the contaminating error in \eqref{DiscreteLinearProb}, the Landweber method can exhibit \emph{semiconvergence} behavior, which manifests visually in the plot of relative reconstruction error norms per iteration. That is, the  error norms decay in early iterations and are followed by a continual monotonic increase. As such, a good termination criterion is needed for the Landweber method. For the numerical experiments, we used a Landweber implementation from the software package AIR Tools II \cite{AIRtools} where the default relaxation parameter $\omega=1.9/\sigma_1^2$ was used.  

Another approach to reduce semiconvergent behavior is to improve the conditioning of the problem by incorporating prior knowledge in the form of regularization. In variational form, the regularized problem may be given by 
\begin{equation} \label{penalizedLS}
\min_\bfx \| \bfA \bfx - \bfb \|^2 + \lambda\,R(\bfx)
\end{equation}
where $R$ is a regularization term and $\lambda>0$ is a regularization parameter that balances the trade-off between the data fit and regularization terms. In standard form Tikhonov regularization, $R(\bf{x})$ is replaced with an $\ell_2$ penalization term so that the penalized LS solution is given by
\begin{equation} \label{l2_l2}
    \bfx_\lambda = \argmin_\bfx \left\| \bfA \bfx - \bfb \right\|^2 + \lambda\left\|\bf{x}\right\|^2.
\end{equation}
For large-scale problems, iterative methods can be used to compute approximations to \eqref{l2_l2}, where at the $k$th iteration, a solution is sought in a $k$ dimensional Krylov subspace.  That is, the $k^{th}$ iterate is given by
\begin{equation} \label{hybridTik}
    \mathbf{x}_\lambda^{(k)} = \argmin_{\mathbf{x}\in \mathcal{R}(\mathbf{V}_k)} \left\| \bfA \bfx - \bfb \right\|^2 + \lambda\left\|\bf{x}\right\|^2,
\end{equation}
where $\mathcal{R}(\cdot)$ denotes the range of a given matrix and $\mathcal{R}(\mathbf{V}_k) = \mathcal{K}_k(\mathbf{A}^T\mathbf{A},\mathbf{A}^T\mathbf{b})$ when the starting guess, $\mathbf{x}^{(0)}$, is the zero vector; see \cite{ChungGazzola_survey,CNO08} for further exposition. We note that $\mathbf{x}^{(k)} \rightarrow \mathbf{x}_{\lambda}$ as $k \rightarrow \min\{m,n\}$ for $\mathbf{A} \in \mathbb{R}^{m \times n}$.  For problems where the regularization parameter $\lambda$ is not known in advance, hybrid projection methods can be used, where the regularization parameter is estimated during the iterative process, resulting in a sequence of projected problems of the form \eqref{hybridTik} with a different regularization parameter at each iteration,  $\lambda^{(k)}$.  See \cite{ChungGazzola_survey} for a survey on hybrid projection methods.

An alternative choice of $R$ in \eqref{penalizedLS} is $\ell_1$ regularization,
\begin{equation} \label{l2_l1}
    \min_\bfx \left\| \bfA \bfx - \bfb \right\|^2 + \lambda\left\|\bf{x}\right\|_1.
\end{equation}
The choice of $\ell_1$ regularization over $\ell_2$ may stem from a desire to promote sparsity in solutions by reducing the sensitivity to outliers. A popular strategy to solve \eqref{l2_l1} is the iterative shrinkage-thresholding algorithm (ISTA) whose general iterative process is given by
\begin{equation} \label{ISTA}
    \mathbf{x}^{(k+1)} = \Psi_{\lambda}\left(\mathbf{x}^{(k)} + 2t\bfA^T\mathbf{r}^{(k)} \right)
\end{equation}
where $t$ is an appropriate step size, $\mathbf{r}^{(k)}=\mathbf{b}-\mathbf{Ax}^{(k)} \in \mathbb{R}^m$ denotes the residual at the $k$th iteration, and $\Psi_{\lambda}: \mathbb{R}^n \to \mathbb{R}^n$ is the shrinkage operator defined by 
\begin{equation}
    \Psi_{\lambda}\left(x_i\right) = (|x_i|-\lambda)_{+}\text{sgn}\,(x_i)
\end{equation}
where $x_i$ denotes the $i^{th}$ entry of $\mathbf{x}$; see \cite{DDM_04} for more details. 
We use the seminal Fast-ISTA (FISTA) method from Beck and Teboulle \cite{FISTA} with backtracking to efficiently compute a solution to \eqref{l2_l1}. FISTA differs from its predecessor through the choice of its shrinkage operator $\Psi_{\lambda}(\cdot)$ which acts on a specific linear combination of the previous two iterates $\left\{\mathbf{x}^{(k-1)},\mathbf{x}^{(k)}\right\}$. This choice results in a superior worst-case complexity of $\mathcal{O}(1/k^2)$ versus $\mathcal{O}(1/k)$ for ISTA. We use the FISTA implementation from the software package IR-Tools \cite{IRtools}.

Obtaining a suitable choice of the regularization parameter ($\lambda$ for variational regularization and the stopping iteration $k$ for Landweber) is a non-trivial task and has continued to receive considerable attention over the last four decades; see references within \cite{Bauer11,Bauer15,RR13}.

Since our focus herein is on the performance of Landweber and other iterative regularization methods to recover cosmic strings, we have used the optimal regularization parameter $\lambda_{opt}$ in each experiment.  That is, given the true solution $\mathbf{x}_{\rm true}$, we determine the optimal regularization parameter to be that which minimizes the relative reconstruction error. For Landweber, this corresponds to finding $k_{opt} = \argmin_k \|\bfx^{(k)} - \bfx_{\rm true}\|$ within a maximal number of iterations. For the solution of \eqref{hybridTik}, we select at iteration $k$
\begin{equation*}
        \lambda^{(k)}_{opt} = \argmin_{\lambda} \frac{\|\mathbf{x}^{(k)}_{\lambda} - \mathbf{x}_{\rm true}\|}{\|\mathbf{x}_{\rm true}\|}.
\end{equation*}
This process can be efficiently implemented at each iteration given that $\mathbf{x}_{\lambda}^{(k)}$ is computed in the subspace $\mathcal{R}(\mathbf{V}_k)$ where typically $k \ll \min\{m,n\}$; see \cite{ChungGazzola_survey, CNO08} for more details. A similar process was carried out involving the minimization of \eqref{l2_l1} to determine the regularization parameter that would result in the lowest relative error amongst the maximal number of iterations within a given implementation of FISTA.

\section{Numerical Results} \label{sec-results}

We illustrate the performance of Landweber, Tikhonov, and FISTA described in Section \ref{sec-iterMeth}  for recovering cosmic strings in a two dimensional universe (i.e.\ $n=2$). To do this, we model a cosmic string as a closed loop in $\mbr^2$ and consider its evolution in $\mbr^{2+1}$. As explained in the introduction, we aim to recover the string through the gravitational perturbations it produces. We take a scalar function $f$ with conormal type singularities to be the string to the world sheet in $\mbr^{2+1}$. From the tomography data $Lf$ we seek to recover $f$. The details of the model are described in Examples $\#1-\#3$ below. We remark that although our model is highly idealized, it does capture some of the main features of the physical problem. We refer the reader to \cite{ViSh} for additional background.

In our experiments we assume knowledge of an upper bound on the error contaminating the problems, i.e.,
\[
\left\|\bf{e}\right\|\leq\delta,
\]
and utilize the discrepancy principle (DP) as a termination criterion for the iterative methods considered. 
Let $\mathbf{x}^{(k)}$ denote the solution determined by the appropriate algorithm at iteration $k$.
The DP prescribes that an iterative method 
should be terminated as soon as an iterate is found that satisfies
\begin{equation}\label{discrp}
\left\|\mathbf{Ax}^{(k)} - \bf{b}\right\| \leq \tau \delta.
\end{equation}
Here, $\tau>1$ is a user specified constant that is independent of 
$\delta$; see \cite{H96,Ha98}. In our experiments we utilize $\tau=1.01$. We iteratively track the relative residual norm (RRN) defined by 
\begin{equation*} 
\text{RRN}\left(\mathbf{x}^{(k)}\right)=\frac{\left\|\mathbf{Ax}^{(k)} - \mathbf{b}\right\|}
{\left\|\bfA\mathbf{x}_{\rm true}\right\|}
\end{equation*}
to determine when the normalized DP is satisfied. To evaluate the quality of the computed solutions, 
we compute the relative reconstructive error (RRE) defined by
\begin{equation*} 
\text{RRE}\left(\mathbf{x}^{(k)}\right)=\frac{\left\|\mathbf{x}^{(k)} - \mathbf{x}_{\rm true}\right\|}
{\left\|\mathbf{x}_{\rm true}\right\|}.
\end{equation*}

The reported numerical results of this work were carried out in MATLAB R2022b 64-bit on a MacBook Pro laptop 
running MacOS Ventura with an Apple M2 Pro processor with @3.49 GHz and 16 GB of RAM. The codes for this paper will be made available at \url{https://github.com/lonisk1/CMB_InvProb} once the revision process is complete.

{\bf Example $\#1$:}  For our first example, we consider a string moving in the positive $x$-direction at constant speed $c\geq 0$. We use $(x, y)$ for coordinates of $\mbr^2$ and $t$ for the time variable. We model the string as the unit circle $\gamma = \{(x, y)\in \mbr^2: x^2 + y^2 = 1\}$ where the world sheet of $\gamma$ is a surface in $\mbr^{2+1}$. The following function with conormal type singularities to the world sheet is given by
\beq
f_1(t, x, y) = (2 - ((x - ct)^2 + y^2)^\ha)_+^a, \quad t\in [0, 4], \quad x,y\in [-3,3].
\eeq 
Here, $a\in \mbr$ controls the strength of conormal singularity, which for this example we fix to be $a=0.05$. We note that the speed of light for our setup is $1$. We consider three cases: (i) $c = 0$; (ii) $c = 1$; (iii) $c=1.5$ corresponding to a string that is stationary, moving at the speed of light, and moving faster than the speed of light. From the analysis in Sections \ref{sec-land} and \ref{sec-art}, we expect to see numerical instabilities in case (ii) and especially in case (iii).

For the discretization of the world sheet, the $x,y$-values are linearly spaced on a $51\times51$ grid, and $t$ is linearly discretized into $20$ slices where the planes are fixed to be the midpoints of the subintervals of $[0,\,4]$ so that the first and last do not coincide with the source and detector planes. The resulting discretized matrix $\mathbf{A}$ from \eqref{DiscreteLinearProb} is $145,880 \times 52,020$. This is the generic setup for Examples $\#1-\#3$. 
To realistically simulate the problem, noise was added to the observations by forming the vector $\mathbf{e}$ with normally distributed random entries with mean zero so that $\mathbf{b} = \mathbf{Ax}_{\text{true}}+\mathbf{e}$; the vector $\mathbf{e}$ is scaled so as to correspond to a specific noise level given by
\begin{equation*}
v = 100\left(\frac{\|\mathbf{e}\|}{\|\mathbf{Ax}_{\text{true}}\|}\right).
\end{equation*}
In this example $v$ is set to be $5\%$. 

\begin{figure}
\begin{minipage}{1\linewidth}
	\centering
	\begin{minipage}{0.28\linewidth}
		\centering
		\includegraphics[width=\linewidth]{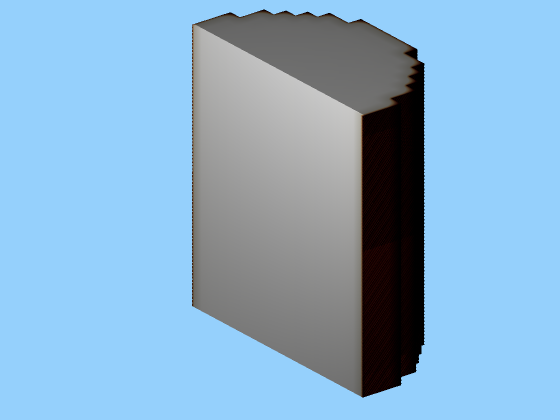}\\$c=0$
	\end{minipage}
	\begin{minipage}{0.28\linewidth}
		\centering
		\includegraphics[width=\linewidth]{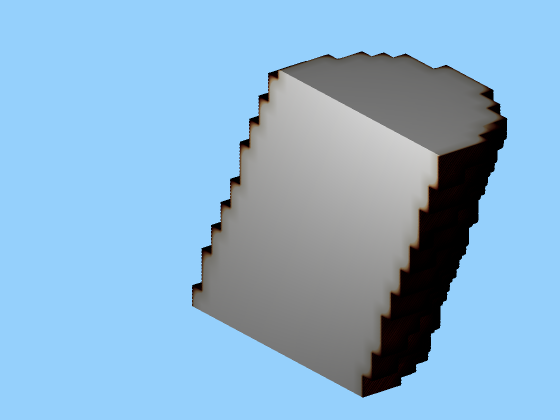}\\$c=1$
	\end{minipage}
	\begin{minipage}{0.28\linewidth}
		\centering
		\includegraphics[width=\linewidth]{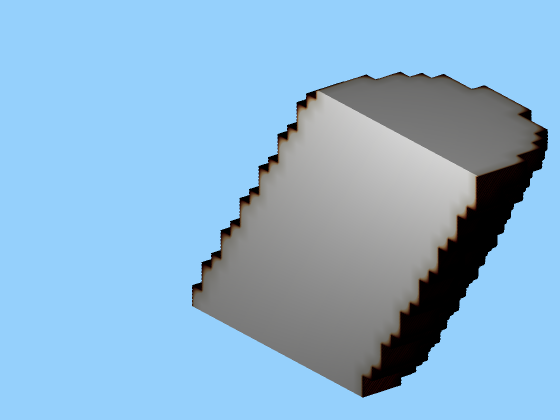}\\$c=1.5$
	\end{minipage}
    \begin{minipage}{0.1\linewidth}
		\centering
		\includegraphics[width=0.85\linewidth]{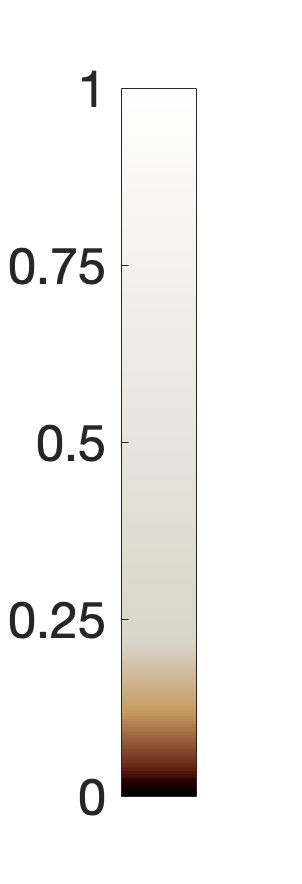}\\ \phantom{$2$} 
	\end{minipage}
\end{minipage}
\caption{Isosurface plots of the true strings of Example $\#$1 with $a=0.05$ for three different values of $c$. A clipping plane through the $yt$-plane is utilized to display the inside of the strings. The temporal direction of each of the strings increases from the bottom of figure to the top.}
\label{trueStrings_ex1}
\end{figure}

The true strings are displayed in Figure \ref{trueStrings_ex1} as an isosurface with a clipping plane cutting halfway through the $yt-\text{plane}$. 
All isosurface plots use an isovalue of zero. The temporal direction of the $3$-dimensional plots in this section increase from the bottom to the top. We point out that the numerical values of the strings are $1$ at the center of any slice and decay towards the edge with the rate of decay controlled by $a$.

\begin{figure}
\begin{minipage}{1\linewidth}
	\centering
	\begin{minipage}{0.495\linewidth}
		\centering
		\includegraphics[width=\linewidth]{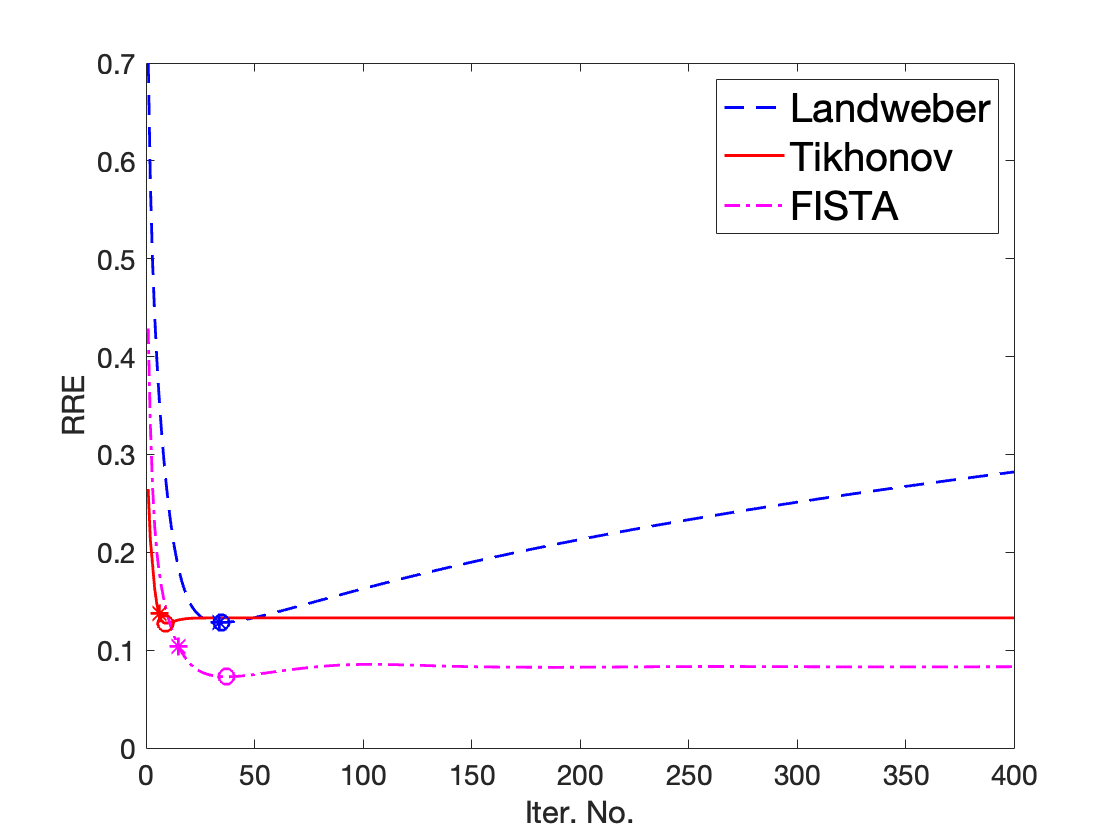}
	\end{minipage}
	\begin{minipage}{0.495\linewidth}
		\centering
		\includegraphics[width=\linewidth]{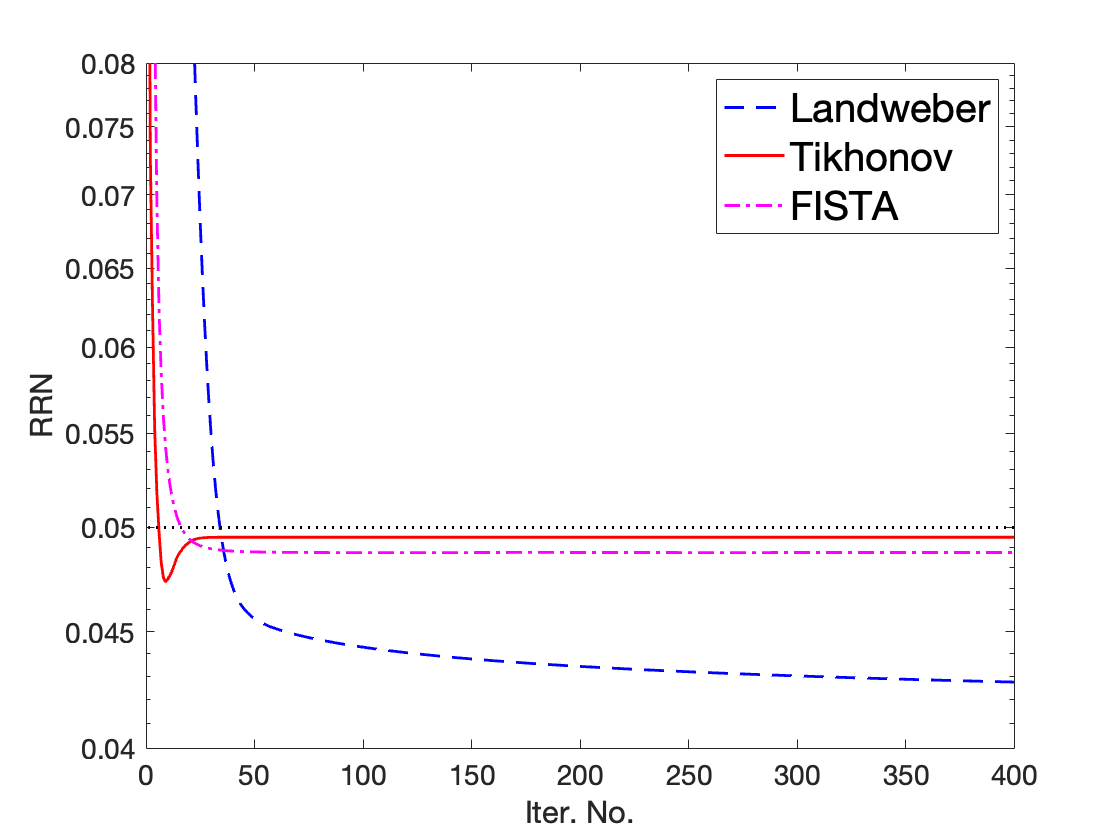}
	\end{minipage}
\end{minipage}
\caption{RRE and RRN plots for Example $\#1$ with $a=0.05$ and $c=0$. The black dashed line in the RRN plot (right) indicates where the methods would terminate according to the DP. Colored stars and circles in the RRE plot (left) indicate termination according to the DP and the smallest RRE achieved in the allotted iterations, respectively.}
\label{plots_ex1}
\end{figure}

The numerical results for each of the cases considered are summarized in Table \ref{table_ex1}. The RRE and RRN plots for the $c=0$ case are shown in Figure \ref{plots_ex1}. The most prominent result is the superior RRE performance of FISTA amongst the considered methods. We note that for all cases except Landweber for $c=0$, utilizing the DP resulted in a suboptimal RRE value, though we stress that knowledge of $\mathbf{x}_{true}$ in real applications is infeasible. It is noteworthy, however, that for Landweber and Tikhonov that the DP RRE values did not deviate significantly from the minimally obtained RRE values within the allotted iterations compared to the relative difference between the RRE value attained by FISTA at the DP and its minimum RRE. As mentioned in Section \ref{sec-iterMeth}, with the optimal choice of the regularization parameter, the RRE curves of Tikhonov and FISTA do not experience the semiconvergence behavior that Landweber does.

\begin{table}
	\centering
		\begin{tabular}{c |c | c | c }
 			c-value & Method & RRE (DP Iter. no.) & Min. RRE (Iter. no.) \\
			\ChangeRT{1.5pt}
			& \footnotesize{Landweber}& $0.1283\,\,(34)$ & $0.1282 \,\,(35)$ \\ 
			$c = 0$ & \footnotesize{FISTA}& $0.1039\,\,(15)$ & $0.0728\,\,(37)$ \\
			& \footnotesize{Tikhonov}& $0.1374\,\,(6)$ & $0.1275 \,\,(9)$ \\
			 \hline
                \hline
			&  \footnotesize{Landweber}& $0.1596\,\,(34)$ & $0.1591 \,\,(39)$ \\
			$c=1$ & \footnotesize{FISTA}& $0.1326\,\,(17)$ & $0.0955 \,\,(54)$ \\
			& \footnotesize{Tikhonov}& $0.1683\,\,(7)$ & $0.1594 \,\,(11)$ \\
			     \hline
    			 \hline
			&  \footnotesize{Landweber}& $0.1969\,\,(35)$ & $0.1950 \,\,(47)$ \\
			$c=1.5$ & \footnotesize{FISTA}& $0.1707\,\,(18)$ & $0.1198 \,\,(83)$ \\
			& \footnotesize{Tikhonov}& $0.2092\,\,(7)$ & $0.1957 \,\,(13)$ \\
			 \hline
		\end{tabular}
\caption{RRE results for Example $\#$1 of a string moving in the positive x-direction at speed $c$ with $5\%$ noise for the three methods considered. RRE values are provided with corresponding iteration number for termination by the DP as well as where the minimal error occurs within the maximal number of iterations allotted (400).}
\label{table_ex1}
\end{table}

\begin{figure}
\begin{minipage}{1\linewidth}
	\centering
	\includegraphics[width=\linewidth]{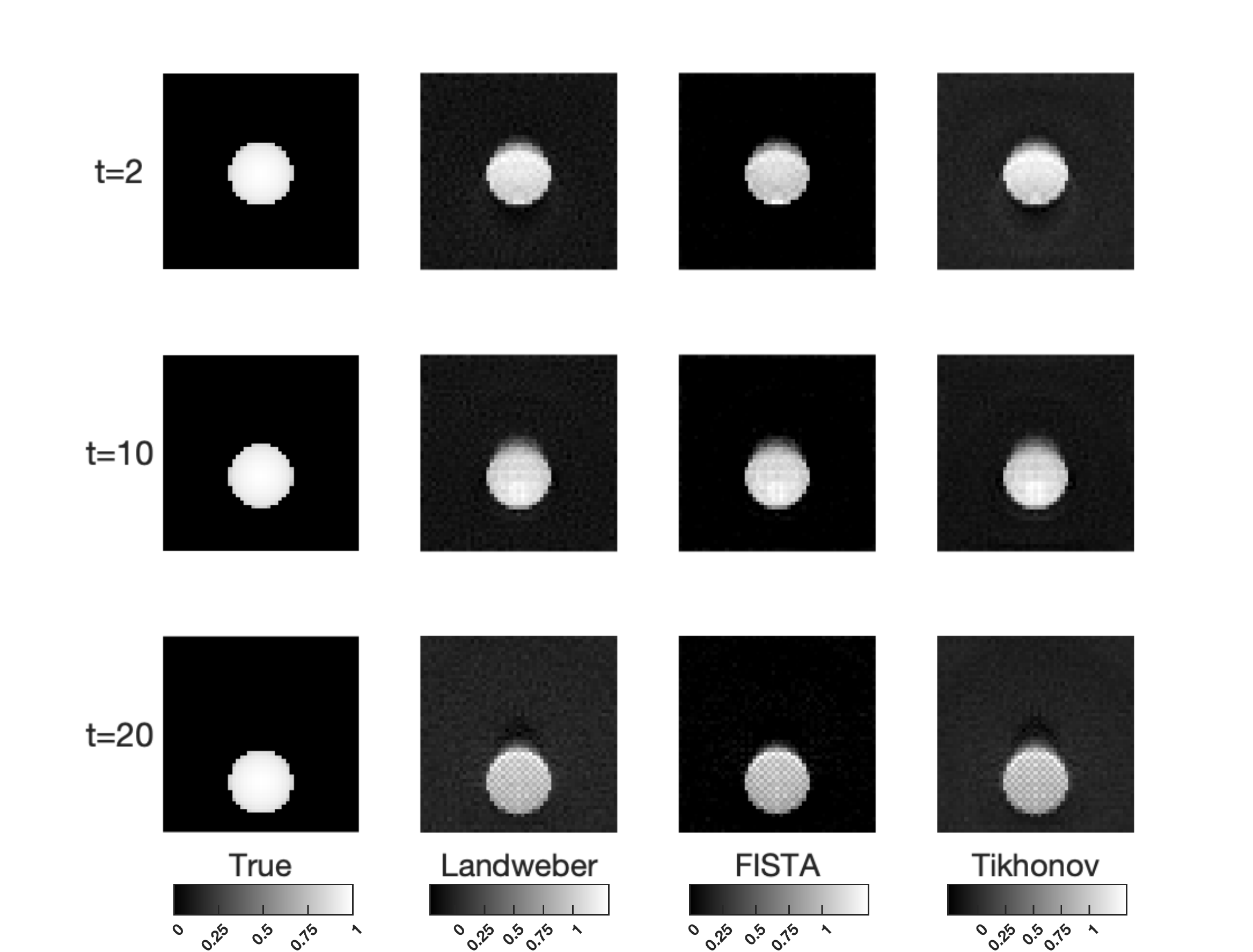}
\end{minipage}
\caption{Reconstructed string for Example $\#1$ when terminated according to the DP at corresponding time slices with $5\%$ noise, $a=0.05$, and $c=1.5$ for each method considered.}
\label{slices_ex1}
\end{figure}

Selected reconstructed slices of the $c=1.5$ case are shown in Figure \ref{slices_ex1} when all methods are terminated according to the DP. It is clear that the parts of the circle near the horizontal directions in the slice are not reconstructed well because the normal directions to these parts of the circle correspond to time-like directions in $\mbr^{2+1}$, which according to the analysis in Section \ref{sec-land} are in the unstable regime for the Landweber method. The same effects also appear in the Tikhonov and FISTA reconstructions. On the other hand, the parts of the circle near the vertical directions in the slices are in the stable regime and the reconstruction is better.

{\bf Example $\#2$:}  
 Let $\gamma$ be the unit circle in Example $\#1$. We consider the scenario where the string collapses to a point at speed $1/2$ and then re-emerges. Consider the following function with conormal singularities to the world sheet,
\beq
f_2(t, x, y) =
\begin{cases}
& (1 -t/2 - (x^2 + y^2)^\ha)_+^a, \quad t\in [0, 2], \quad x,y\in [-3, 3]\\
& (1 -t/2 + (x^2 + y^2)^\ha)_+^a, \quad t\in [2, 4], \quad x,y\in [-3, 3]
\end{cases}
\eeq
Here, we fix $a=1$ so the function is more regular than $f_1$. This is because near the collapsing point, the string is harder to capture on the same grid. The discretization of the world sheet and the addition of $5\%$ noise to the problem is the same as in Example $\#1$.

\begin{figure}
\centering
\begin{minipage}{0.85\linewidth}
	\centering
	\begin{minipage}{0.45\linewidth}
		\centering
		\includegraphics[width=\linewidth]{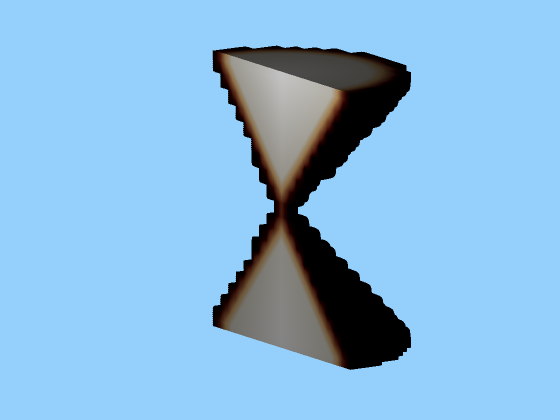}\\True \\ \phantom{} 
	\end{minipage}
	\begin{minipage}{0.45\linewidth}
		\centering
		\includegraphics[width=\linewidth]{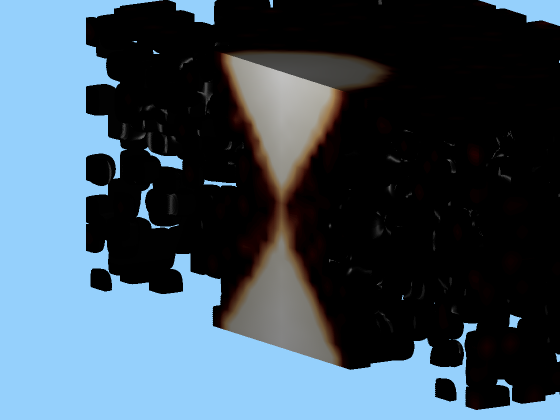}\\FISTA \\ \phantom{} 
	\end{minipage}
    \begin{minipage}{0.45\linewidth}
		\centering
		\includegraphics[width=0.5\linewidth]{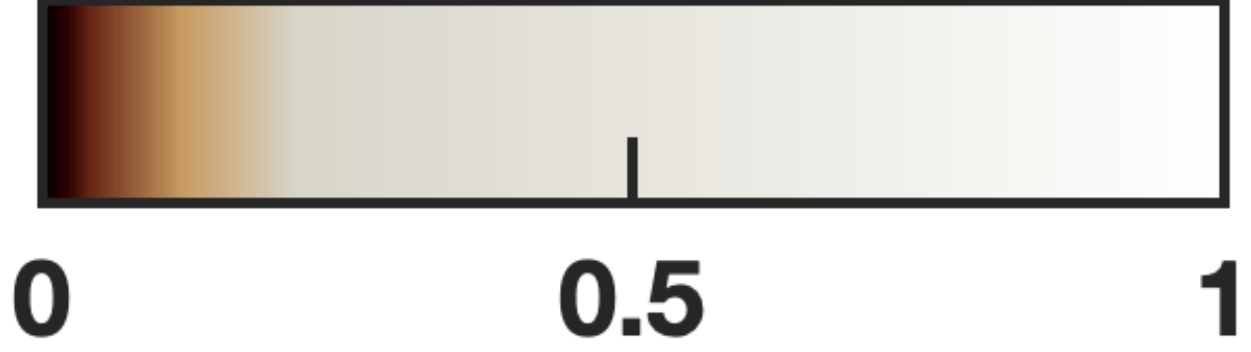}\\ \phantom{} 
	\end{minipage}
    \begin{minipage}{0.45\linewidth}
		\centering
		\includegraphics[width=0.5\linewidth]{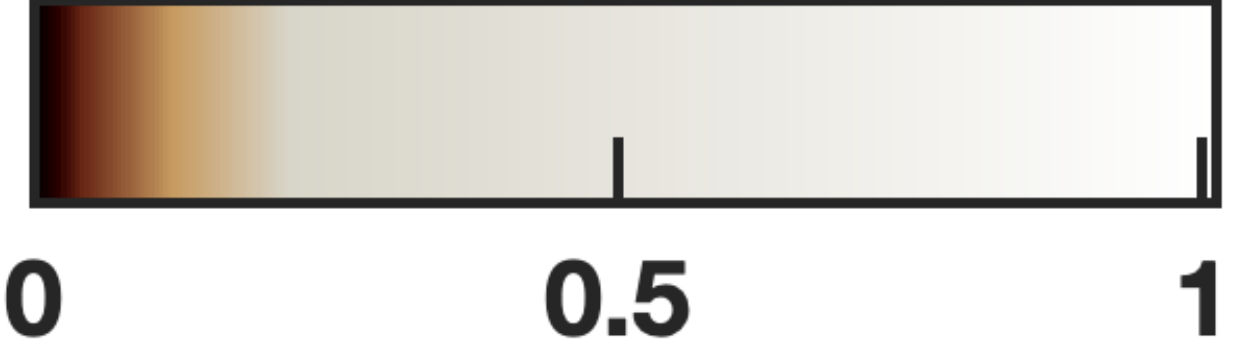}\\ \phantom{} 
	\end{minipage}
    \begin{minipage}{0.45\linewidth}
		\centering
		\includegraphics[width=\linewidth]{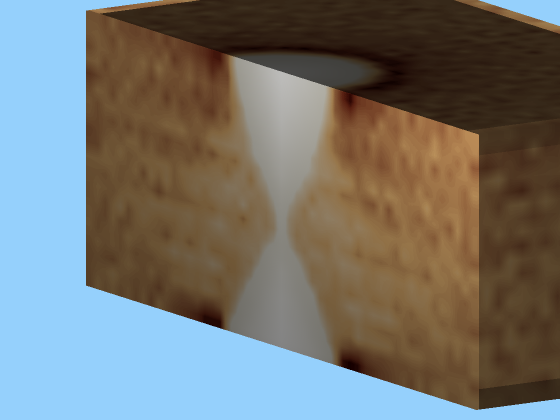}\\Landweber \\ \phantom{} 
	\end{minipage}
	\begin{minipage}{0.45\linewidth}
		\centering
		\includegraphics[width=\linewidth]{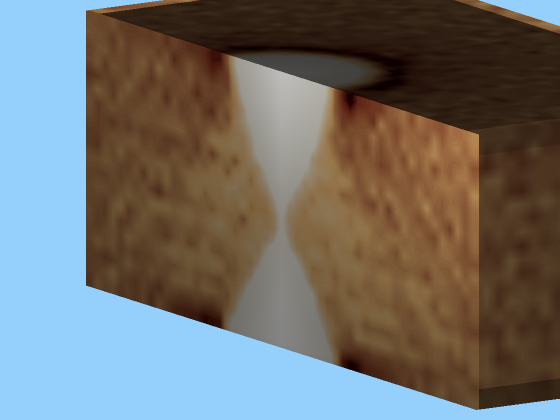}\\Tikhonov \\ \phantom{} 
	\end{minipage}
    \begin{minipage}{0.45\linewidth}
		\centering
		\includegraphics[width=0.5\linewidth]{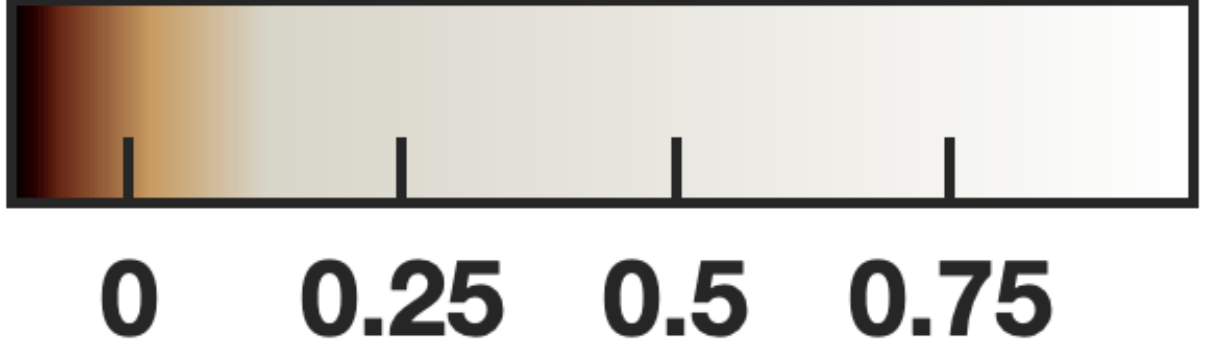}\\ \phantom{} 
	\end{minipage}
    \begin{minipage}{0.45\linewidth}
		\centering
		\includegraphics[width=0.5\linewidth]{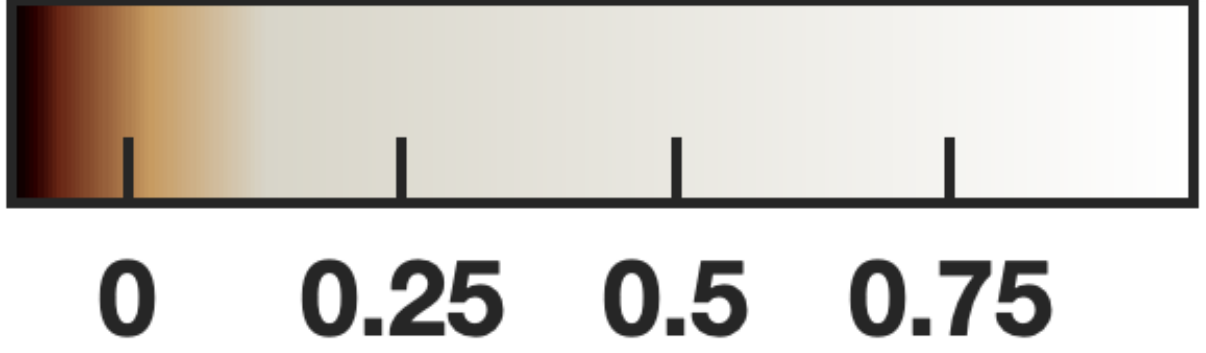}\\ \phantom{} 
	\end{minipage}
 \end{minipage}
\caption{String reconstructions corresponding to lowest RRE for Example $\#$2 via the considered methods for $5\%$ noise with $a=1$. The strings are sliced through the $yt$-plane to provide an interior viewpoint with the temporal direction of each of the strings increases from the bottom of figure to the top.}
\label{3dImages_ex2}
\end{figure}

\begin{figure}
\begin{minipage}{1\linewidth}
	\centering
	\includegraphics[width=\linewidth]{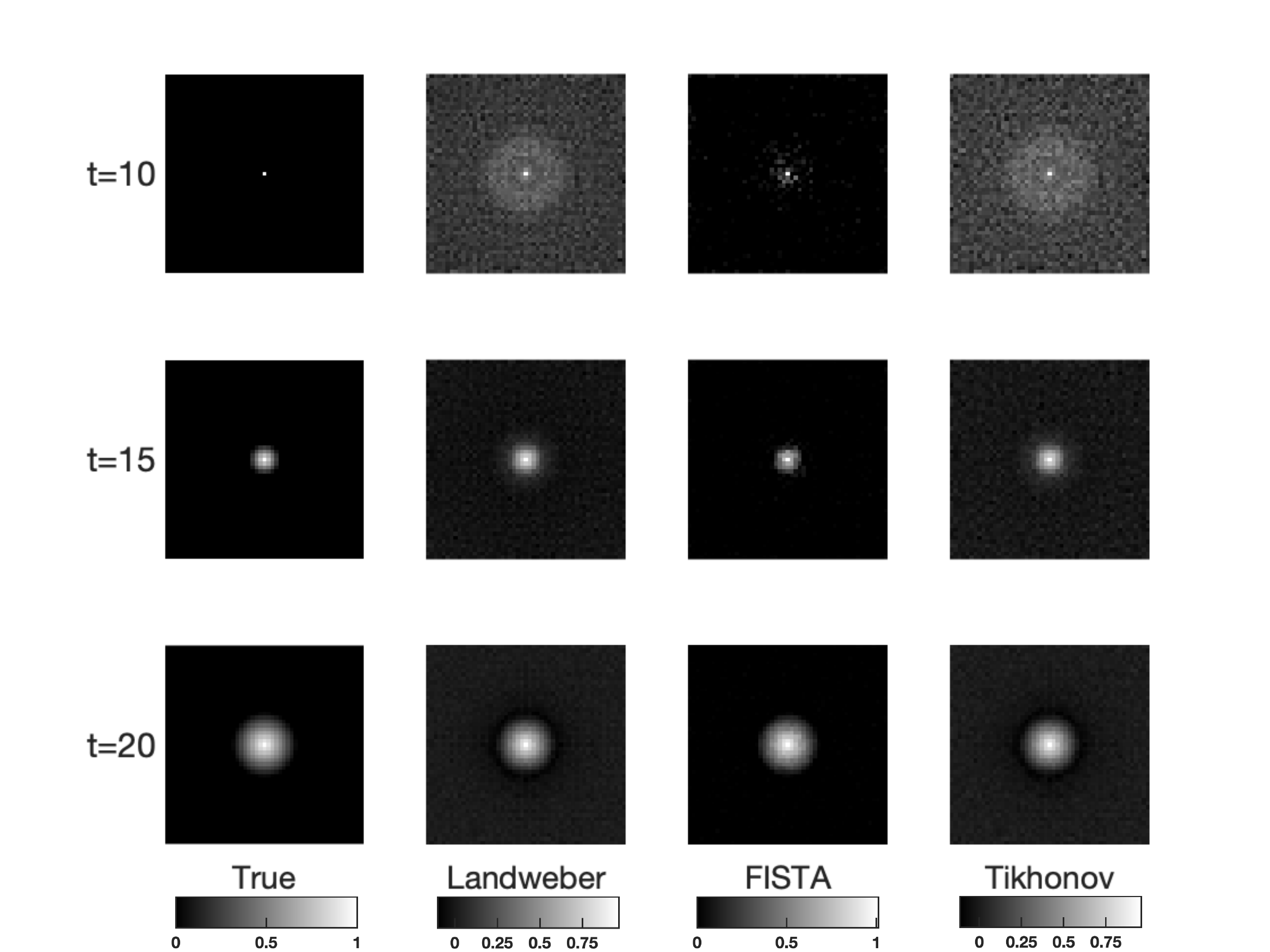}
\end{minipage}
\caption{Selected slices of the lowest RRE string reconstructions for Example $\#2$ with $5\%$ noise and $a=1$.}
\label{slices_ex2}
\end{figure}

The true string and the best reconstructions in terms of lowest RRE values attained are shown in Figure \ref{3dImages_ex2}. Identically to Figure \ref{trueStrings_ex1} from Example $\#1$, the isosurface plots contain a clipping plane cutting halfway through the $yt$-plane. Selected slices are shown in Figure \ref{slices_ex2} for the lowest RRE reconstructions of each method considered. From these two figures, we observe that image reconstruction of the string by FISTA most closely matches the true image of the string. 
From the analysis in Section \ref{sec-land} and \ref{sec-art}, we expect instability and artefacts at the collapsing point. It seems that the singularity of $f_2$ is not strong enough to show the effects. However, near the collapsing point, there are increasing discretization errors which produce artefacts in the light ray directions. These could explain the disk-like artefacts in  the first row of Figure \ref{slices_ex2} for the Landweber and Tikhonov reconstructions.

\begin{table}
	\centering
		\begin{tabular}{c | c | c }
 			Method & RRE (DP Iter. no.) & Min. RRE (Iter. no.) \\
			\ChangeRT{1.5pt}
			\footnotesize{Landweber}& $0.4783\,\,(93)$ & $0.4147 \,\,(400)$ \\ 
			\footnotesize{FISTA}& $0.2566\,\,(43)$ & $0.0988\,\,(226)$ \\
			\footnotesize{Tikhonov}& $0.4754\,\,(15)$ & $0.4043 \,\,(103)$ \\
			 \hline
		\end{tabular}
\caption{RRE results for Example $\#2$ of a string collapsing and re-emerging with $5\%$ noise for the three methods considered. RRE values are provided with corresponding iteration number for termination by the DP as well as where the minimal error occurs within the maximal number of iterations allotted (400).}
\label{table_ex2}
\end{table}

The RRE values for Example $\#2$ are summarized in Table \ref{table_ex2}. Here, unlike for Example $\#1$, the lowest RRE values attainable within the allotted number of iterations for all methods is significantly lower than the corresponding RRE values attained when terminating the methods via the DP; because of this, we chose to highlight the lowest RRE reconstructions amongst the allotted iterations considered in Figures \ref{3dImages_ex2} and \ref{slices_ex2}. For Landweber and Tikhonov, the loss in relative accuracy is approximately $6\%$, and is $15\%$ for FISTA. These accuracy results compared with Example $\#1$ suggest that the DP many not be an appropriate algorithmic termination criterion when the string involved changes significantly with respect to time.

{\bf Example $\#3$:}
For our final example, we consider a string that appears at $t_0$ and disappears at $t_1$. For $x,y\in [-3,3]$, we take 
\beq
f_3(t, x, y) = \begin{cases}
&(1.5 - (x^2 + y^2 +(t-2)^2)^\ha)_+^a, \quad t\in [t_0, t_1]\\
&0, \quad t\in [0, t_0) \cup (t_1, 4]
\end{cases}
\eeq
As in Example $\#1$, we take $a=0.05$.  The main goal of the example is to investigate the reconstruction at the temporal interfaces $t = t_0$ and $t= t_1$ where $f_3$ is singular in time-like directions. This should be quite challenging according to the analysis in Section \ref{sec-land}, but we will see that certain regularizers can reasonably capture the temporal singularity.

The discretization of the world sheet is the same as in the previous two examples. Since the string only exists from $t \in [t_0,\,t_1]$, this corresponds to the top and bottom $6$ temporal slices being empty. On $[0,\,4]$ with $20$ slices this corresponds to $t_0 = 1.2$ and $t_1 = 2.8$. The problem is considered with $5\%$ added noise in the same manner as was done in the previous examples.

Since this example involves slices that we wish to avoid contaminating during reconstruction, it is natural to consider a regularizer $R$ from \eqref{penalizedLS} which penalizes interaction between sequential time slices, $\mathbf{X}(t)$, for $t=1,\dots,T$. Such a regularizer was discussed in \cite{AKS_Chung_genHyBR} for Bayesian inverse problems. The penalized LS problem may be written as
\begin{equation} \label{genHyBR}
    \min_\bfx \| \bfA \bfx - \bfb \|^2 + \lambda \|\bfx \|_{\bfQ^{-1}}^2
\end{equation}
where $\|\bfx \|_\bfM^2 = \bfx\t \bfM \bfx$ for any symmetric positive definite matrix $\bfM.$ 

We note that for dynamic inverse problems, the prior covariance matrix $\bfQ\in \bbR^{n^2T \times n^2T}$ may be very large and dense, and working with its inverse or square root can be computationally difficult.  Generalized Golub-Kahan based iterative methods were described in \cite{AKS_Chung_genHyBR,chung2018efficient}, where the $k^{th}$ iterate is given by 
\[
\bfx^{(k)} = \bfQ \bfy_k \quad \mbox{where} \quad
\bfy_k = \argmin_{\bfy \in \mathcal{R}(\bfV_k)} \| \bfA \bfQ \bfy - \bfb \|^2 + \lambda \|\bfy\|_\bfQ^2
\]
where $\mathcal{R}(\bfV_k) = \calK_k(\bfA\t \bfA \bfQ, \bfA\t \bfb)$. 

For the numerical experiments, we let $\bfQ$ represent a three-dimensional Mat{\' e}rn covariance kernel with $\nu = 1.5$ and $\ell = 0.05$. This is a nonseparable covariance kernel that models space-time interactions of variability.  Constructing such a matrix is not done explicitly; instead, circulant embedding and FFT-based techniques are used for efficient matrix-vector products. We will refer experimentally to this method as generalized Tikhonov (gen. Tikhonov) where like past examples, we will consider the optimal regularization parameter at each iteration.

\begin{table}
	\centering
		\begin{tabular}{c | c | c }
 			Method & RRE (DP Iter. no.) & Min. RRE (Iter. no.) \\
			\ChangeRT{1.5pt}
			\footnotesize{Landweber}& $0.6518\,\,(44)$ & $0.6300 \,\,(397)$ \\ 
			\footnotesize{FISTA}& $0.3460\,\,(59)$ & $0.2875\,\,(140)$ \\
			\footnotesize{Tikhonov}& $0.6521\,\,(10)$ & $0.6256 \,\,(400)$ \\
            \footnotesize{gen. Tikhonov}& $0.4586\,\,(180)$ & $0.4513 \,\,(317)$ \\
			 \hline
		\end{tabular}
\caption{RRE results for Example $\#3$ of a string with speed $c=0$ with $5\%$ noise for the methods considered. RRE values are provided with corresponding iteration number for termination by the DP as well as where the minimal error occurs within the maximal number of iterations allotted (400).}
\label{table_ex3}
\end{table}

\begin{figure}
\centering
	\begin{minipage}{0.55\linewidth}
		\centering
		\includegraphics[width=\linewidth]{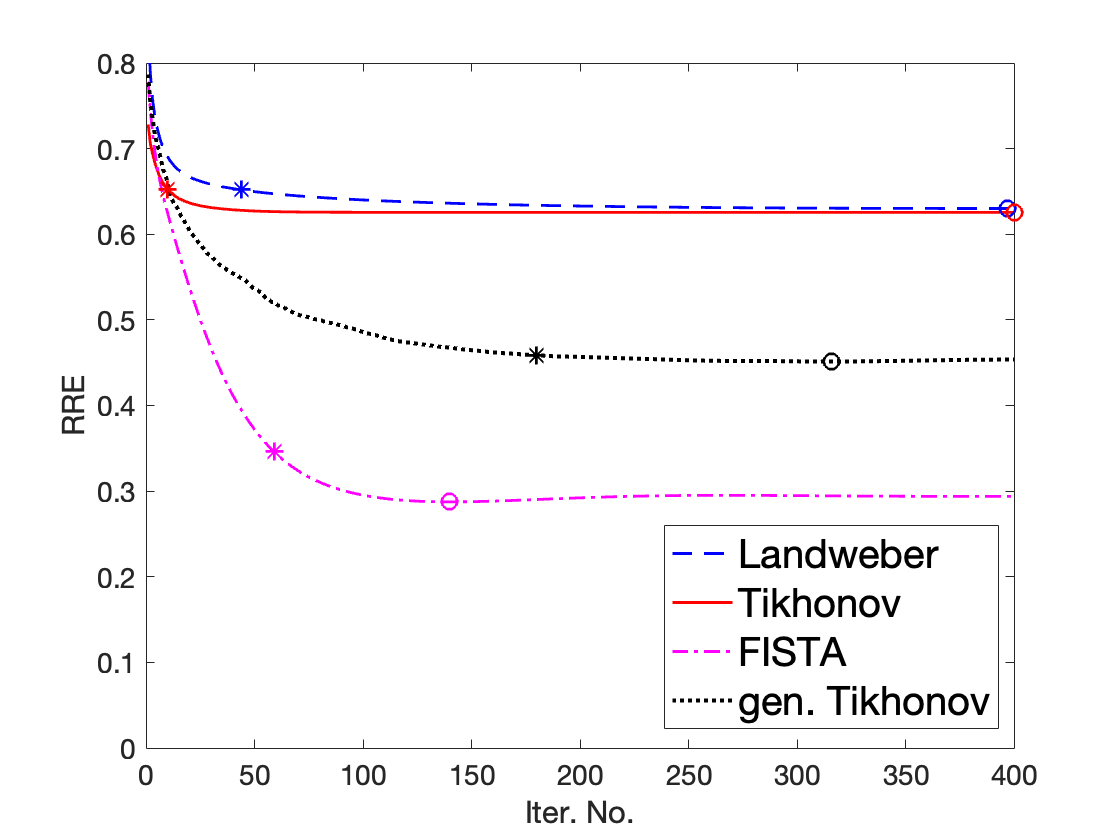}
	\end{minipage}
\caption{RRE plot for Example $\#3$ with $a=0.05$. Colored stars and circles indicate termination according to the DP and the smallest RRE achieved in the allotted iterations, respectively.}
\label{plot_ex3}
\end{figure}

The RRE values for Example $\#3$ provided in Table \ref{table_ex3} provide evidence that the FISTA method provides superior reconstructive power given sufficient iterations with the difference between its lowest RRE and RRE achieved by the DP to be approximately $6\%$. Similarly to Example $\#1$, we do not see significant improvement of the optimal RRE values for Landweber and Tikhonov over their values when terminated according to the DP. As was expected, the gen. Tikhonov method with the Matern covariance regularizer resulted in numerically lower RRE values compared with standard Tikhonov. The RRE values per iteration for all methods are provided in Figure \ref{plot_ex3}. 
We remark that $\bfQ$ and its hyperparameters were empirically chosen and set; further work is needed to determine a good prior covariance for these problems.

\begin{figure}
\centering
\begin{minipage}{0.85\linewidth}
	\centering
	\begin{minipage}{0.45\linewidth}
		\centering
		\includegraphics[width=\linewidth]{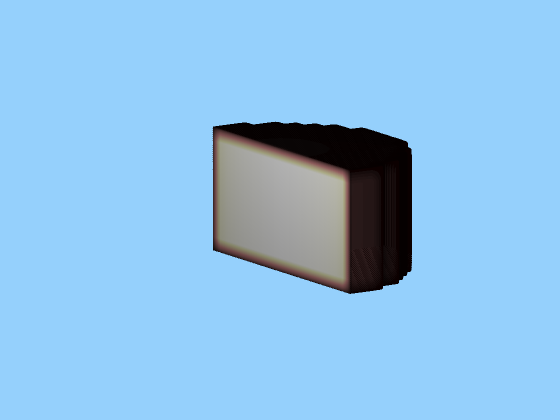}\\True \\ \phantom{} 
	\end{minipage}
	\begin{minipage}{0.45\linewidth}
		\centering
		\includegraphics[width=\linewidth]{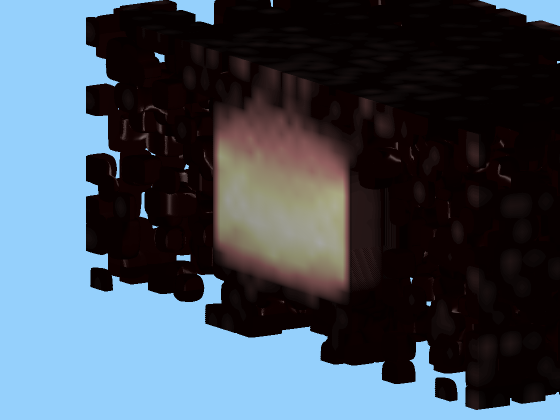}\\FISTA \\ \phantom{} 
	\end{minipage}
    \begin{minipage}{0.45\linewidth}
		\centering
		\includegraphics[width=0.5\linewidth]{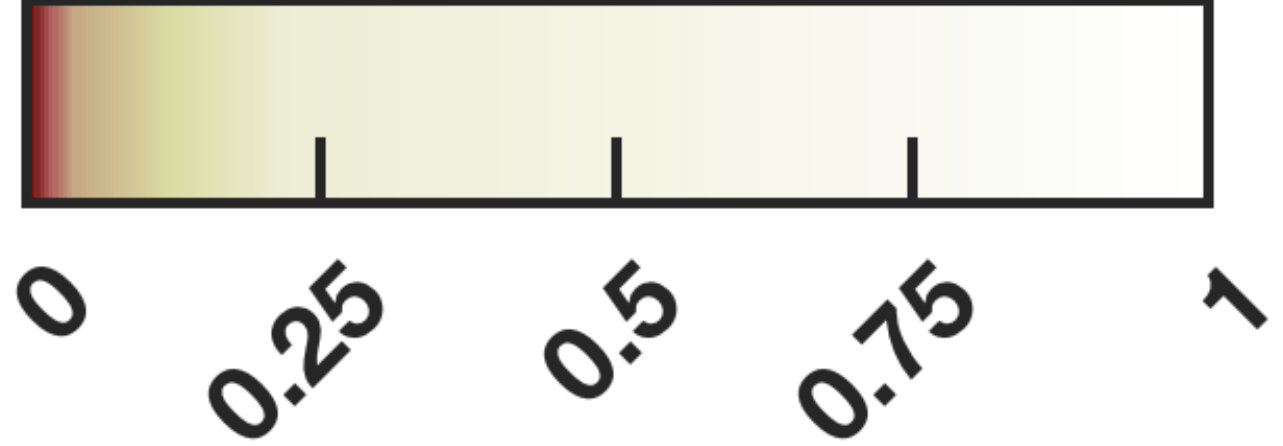}\\ \phantom{} 
	\end{minipage}
    \begin{minipage}{0.45\linewidth}
		\centering
		\includegraphics[width=0.5\linewidth]{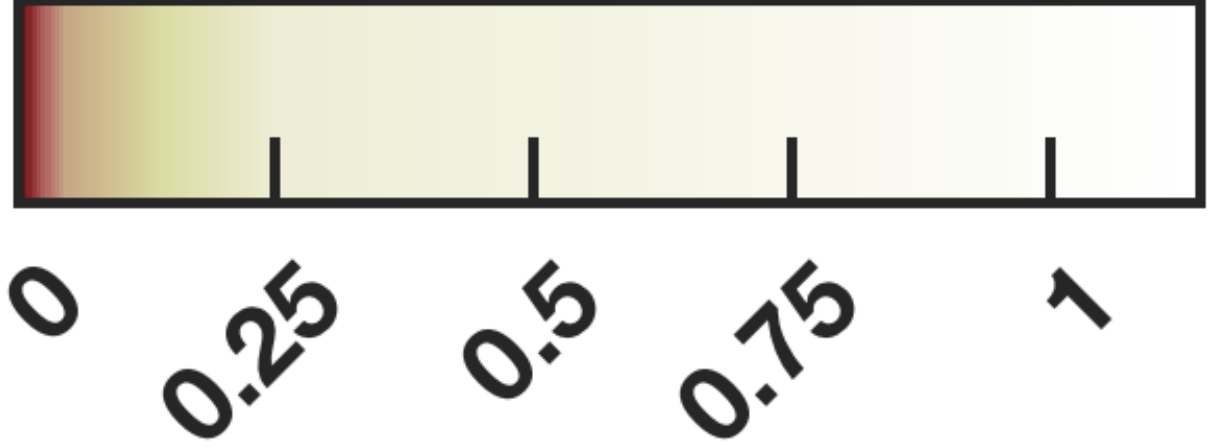}\\ \phantom{} 
	\end{minipage}
    \begin{minipage}{0.45\linewidth}
		\centering
		\includegraphics[width=\linewidth]{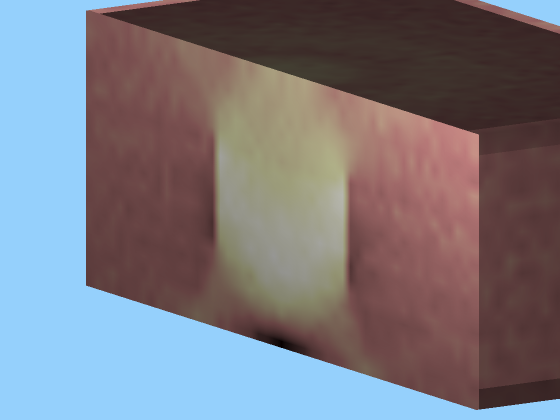}\\Landweber \\ \phantom{} 
	\end{minipage}
	\begin{minipage}{0.45\linewidth}
		\centering
		\includegraphics[width=\linewidth]{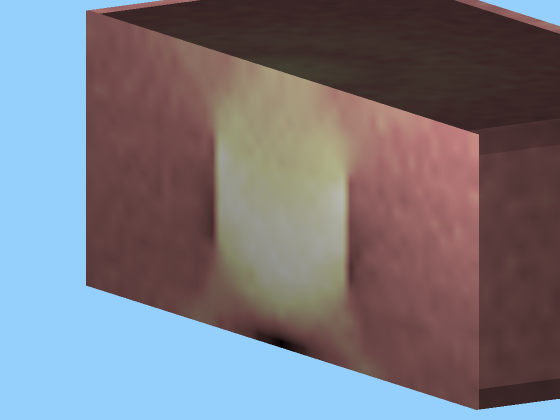}\\Tikhonov \\ \phantom{} 
	\end{minipage}
    \begin{minipage}{0.45\linewidth}
		\centering
		\includegraphics[width=0.5\linewidth]{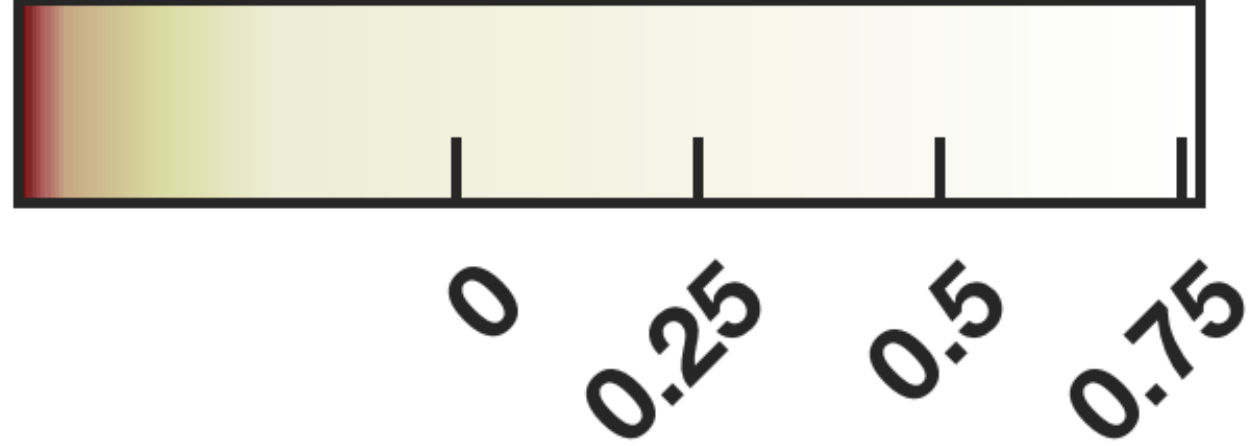}\\ \phantom{} 
	\end{minipage}
    \begin{minipage}{0.45\linewidth}
		\centering
		\includegraphics[width=0.5\linewidth]{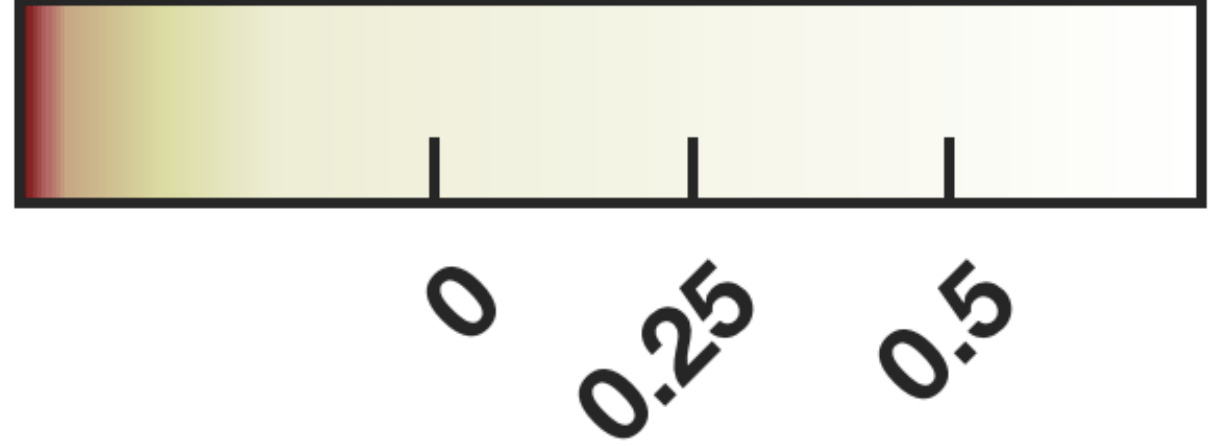}\\ \phantom{} 
	\end{minipage}
 \begin{minipage}{0.45\linewidth}
		\centering
		\includegraphics[width=\linewidth]{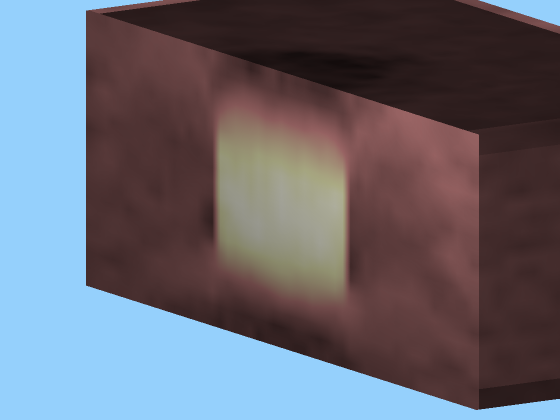}\\gen. Tikhonov \\ \phantom{}
	\end{minipage}
 \begin{minipage}{0.45\linewidth}
		\centering
        \phantom{}
	\end{minipage}
 \begin{minipage}{0.45\linewidth}
		\centering
		\includegraphics[width=0.5\linewidth]{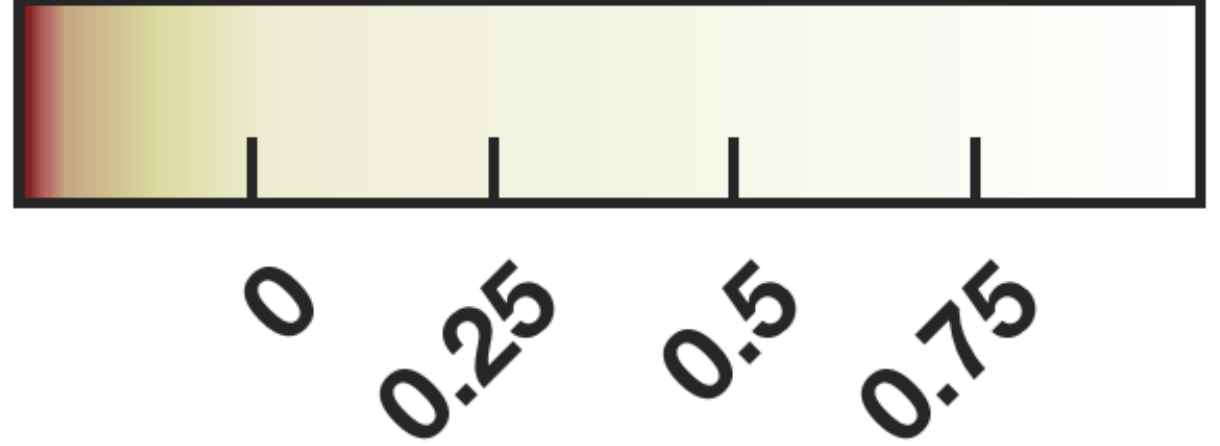}\\ \phantom{}
	\end{minipage}
 \begin{minipage}{0.45\linewidth}
		\centering
        \phantom{}
	\end{minipage}
 \end{minipage}
\caption{String reconstructions corresponding to lowest RRE for Example $\#$3 via the considered methods for $5\%$ noise with $a=0.05$ and $c=0$. The strings are sliced through the $yt$-plane to provide an interior viewpoint with the temporal direction of each of the strings increasing from the bottom of the figure to the top.}
\label{3dImages_ex3}
\end{figure}

\begin{figure}
\begin{minipage}{1\linewidth}
	\centering
	\includegraphics[width=\linewidth]{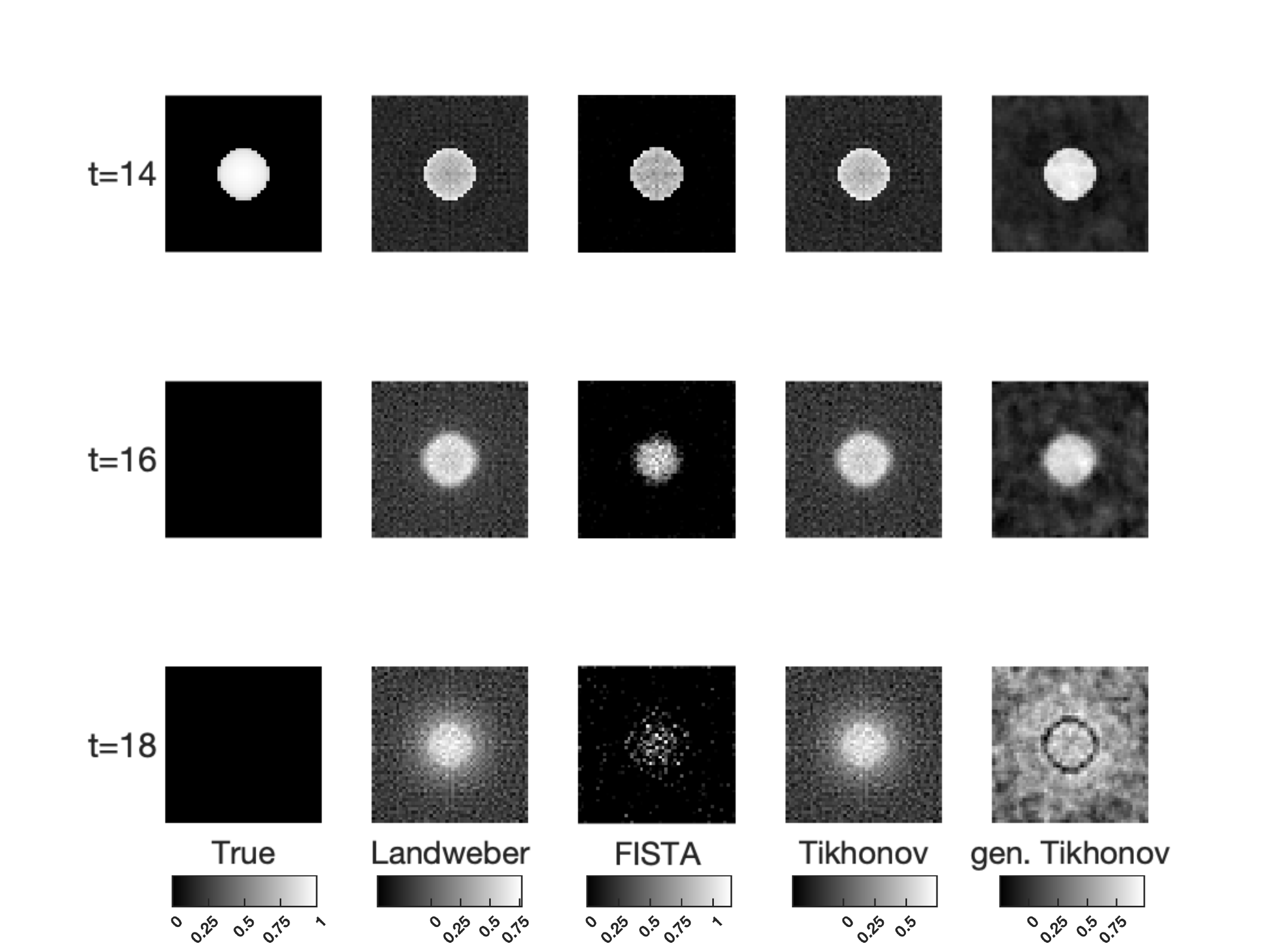}
\end{minipage}
\caption{Selected slices of the lowest RRE string reconstructions for Example $\#3$ with $5\%$ noise and $a=0.05$ and $c=0$ for each method considered.}
\label{slices_ex3}
\end{figure}

The true string and the volumetric reconstructions corresponding to the optimal RRE for each method are shown in Figure \ref{3dImages_ex3} with a clipping plane cutting halfway through the $yt$-plane.  As predicted by the analysis in Section \ref{sec-land}, the temporal boundaries $t = t_0$ and $t=t_1$ are difficult to reconstruct using the Landweber and standard Tikhonov methods. Moreover, there are visible artefacts from the singular support of $f$ at $t = t_0, t_1$ in the light ray directions, which agrees with the analysis in Section \ref{sec-art}. We point out the visible conic structure in the Landweber and Tikhonov reconstructions as was expected from the analysis in Section \ref{sec-art} (see Figure \ref{fig-sphere}). The FISTA and gen. Tikhonov methods performed the best as their reconstructions contain significantly fewer artefacts in the first and last $6$ slices. Selected reconstructed slices from the optimal reconstruction from each method are shown in Figure \ref{slices_ex3}. Here, we chose to present the $16$th and $18$th slices to highlight the existence of reconstruction artefacts from the Landweber and Tikhonnov methods, which are less prominent in the reconstructions of FISTA and gen. Tikhonov. We also provide the $14$th slice to emphasize the transition phase of the string.

\section{Conclusion}
\label{sec-conc}
In this paper, we provide a theoretical and numerical study of iterative methods for the inverse problem of the light ray transform, and discuss applications in cosmology. Specifically, we investigated the use of state-of-the-art iterative regularization methods. We provided a theoretical analysis of the Landweber method and its expected behavior for the recovery of cosmic strings. We numerically compared Landweber against Tikhonov and $\ell_1$ iterative regularization methods. Various examples in Section \ref{sec-results} complement the analysis from Section \ref{sec-art}, demonstrating that the recovery of functions with conormal type singularities (e.g., strings to the world sheet) can be difficult and can lead to reconstruction artefacts. 

Future work includes extensions to a three-dimensional universe over time (i.e., 4D reconstructions). This will require additional tools (e.g., structure exploiting methods and algorithms that can handle larger datasets).
Moreover, for a more realistic setting, the CMB data is measured near a freely falling observer, see for example \cite{LOSU}. We plan to investigate the partial data problem for the light ray transform in the near future. The problem is conceivably more challenging. However, recent studies in \cite{VaWa, Wan1} have shown that stability of the inverse problem can be improved by taking into account the evolution of the universe. It would be interesting to investigate the effects on the numerical reconstruction.
 
\section*{Acknowledgments} 
This work was partially supported by the National Science Foundation (NSF) under grants DMS-2341843 (J. Chung), DMS-2038118 (L. Onisk), and DMS-2205266 (Y. Wang).

\bibliographystyle{siam}
\bibliography{references}

\end{document}